\documentclass[12pt]{amsart}
\usepackage{tikz-cd} 
\usepackage[left=25mm,right=25mm, top=25mm, bottom=25mm]{geometry}
\usepackage[english]{babel}
\usepackage{dirtytalk}
\usepackage[utf8x]{inputenc}
\usepackage[T1]{fontenc}
\usepackage{amsmath}
\usepackage{amstext}
\usepackage{amsfonts}
\usepackage{amssymb}
\usepackage{amsthm}
\usepackage{mathtools}
\usepackage{dsfont}
\usepackage{color}
\usepackage{graphicx}
\usepackage[colorinlistoftodos]{todonotes}
\usepackage[colorlinks=true, allcolors=blue]{hyperref}
\usepackage{float}
\usepackage[colorinlistoftodos]{todonotes}
\usepackage{theoremref}
\usepackage{tikz-cd}
\usepackage[all]{xy}
\allowdisplaybreaks
\let\oldaddcontentsline\addcontentsline

\newcommand{\starttocentries}{\let\addcontentsline\oldaddcontentsline}
\newtheorem{theorem}{Theorem}[section]
\newtheorem{lemma}[theorem]{Lemma}
\newtheorem{proposition}[theorem]{Proposition}
\newtheorem{corollary}[theorem]{Corollary}

\newtheorem*{cor*}{Corollary}
\newtheorem*{thm*}{Theorem}
\newtheorem*{lem*}{Lemma}

\newtheorem*{prop*}{Proposition}
\theoremstyle{definition}
\newtheorem{definition}[theorem]{Definition}

\newtheorem{example}[theorem]{Example}
\newtheorem*{defn*}{Definition}

\theoremstyle{remark}
\newtheorem{remark}[theorem]{Remark}

\newcommand{\G}{G}                           
\newcommand{\g}{g}

\newcommand{\bnd}{\operatorname{{ bnd}}}

\newcommand{\supp}{\operatorname{supp}}
\newcommand{\id}{\operatorname{id}}

\newcommand{\Sub}{\operatorname{Sub}}

\newcommand{\prob}{\text{Prob}}
\newcommand{\poi}{\text{Poi}}

\newcommand{\pr}{\text{pr}}

\newcommand{\bary}{\text{bar}}

\title[Relative stationary systems]{Relative stationary dynamical systems}

\author[T. Amrutam]{Tattwamasi Amrutam} %\orcid{0000-0002-0691-5821}
\address{T. Amrutam \newline IMPAN, 
 00-656 Warszawa, Poland
}
\email{tattwamasiamrutam@gmail.com}

\author[M. Kl\"otzer]{Martin Kl\"otzer}
\address{M. Kl\"otzer \newline Universität Innsbruck,
	Department of Mathematics.
	6020 Innsbruck, Austria.
}
\email{martin.kloetzer@uibk.ac.at}

\author[H. Oppelmayer]{Hanna Oppelmayer} %\orcid{0000-0002-1998-5407} 
\address{H. Oppelmayer \newline Universität Innsbruck,
	Department of Mathematics.
	6020 Innsbruck, Austria.
}
\email{hanna.oppelmayer@uibk.ac.at}

\subjclass[2020]{Primary 37A50; Secondary , 22F10, 37B05, 60J50}
\keywords{Random walks on groups; Poisson boundary; Stationary Dynamics}

\date{\today}

\begin{document}

\begin{abstract} 
Let $G$ be a locally compact second countable group equipped with an admissible non-degenerate Borel probability measure $\mu$.
We generalise the notion of $\mu$-stationary systems to $\mu$-stationary $\G$-factor maps $\pi: (X,\nu)\to (Y,\eta)$.
For these stationary relations between dynamical systems, we provide a structure theorem, which generalises the structure theorem of Furstenberg-Glasner. Furthermore, we show the existence and uniqueness of a relative version of the Poisson boundary in this setup.
\end{abstract}
\maketitle
\tableofcontents

\section{Introduction}

The study of random walks on groups and their associated harmonic functions has been a cornerstone of modern dynamical systems, linking probability theory, group theory, geometry, and operator algebras (see e.g. \cite{FurstenbergPoiFormula}, \cite{FurstenbergDiscrete}, \cite{FurstenbergLie}, \cite{KaimPoi}, \cite{KaimVershikDiscr}, \cite{WoessBoundaryEnds}, \cite{HartKal}). At the heart of this theory lies the notion of a \textit{$\mu$-stationary measure}. Given a locally compact second countable group $\G$ and an admissible probability measure $\mu$, a $\G$-space $(X,\nu)$ is $\mu$-stationary if the measure $\nu$ is invariant under the $\mu$-averaged transition operator, i.e.
$$\nu=\int_G g\nu\, d\mu(g).
$$ This generalises the concept of a measure-preserving action, providing a framework where invariant measures may not exist — such as actions on boundaries of hyperbolic groups or symmetric spaces — but stationary measures are guaranteed to exist in compact setups.

Foundational work by Furstenberg introduced the \textit{Poisson boundary}, a probability space that captures the asymptotic behaviour of random walks on $\G$ (see \cite{FurstenbergNoncomm}, \cite{FurstenbergPoiFormula}). The Poisson boundary is the maximal $\mu$-proximal $\mu$-stationary space, effectively encoding all bounded $\mu$-harmonic functions on the group. This object has proven indispensable in rigidity theory, the study of lattice subgroups, and the structure of $C^*$-algebras associated with groups.

While the absolute theory of stationary systems is well-developed, deep questions remain regarding the \textit{relative} setting. In ergodic theory, the study of factor maps $\pi: (X,\nu) \to (Y,\eta)$ (or extensions) often yields structural insights that are invisible when studying single systems in isolation. The relative theory for measure-preserving actions was notably advanced %by Furstenberg \cite{FurstenbergHomog} and later 
by Zimmer \cite{ZimmerInduced} and Nevo-Zimmer \cite{NZhomog} in the context of semisimple Lie groups as a key tool to study induced actions and rigidity in higher rank setups. However, the category of stationary actions has lacked a parallel relative structural theory so far.

In this paper, we introduce and systematically analyse the concept of \textit{relative $\mu$-stationarity}. A factor map $\pi:(X,\nu)\longrightarrow (Y,\eta)$ between non-singular $G$-spaces is relatively $\mu$-stationary if the fiber measures (disintegrations) $\{\nu_y\}_{y\in Y}$ satisfy invariance on average:
$$\nu_y=
\int_{\G} \g\nu_{\g^{-1}y} \, d\mu(\g)
$$
for $\eta$-almost every $y\in Y$. This definition is a natural interpolation: if $Y$ is a point, we recover classical $\mu$-stationarity; if the equation holds for every $g \in G$ (pointwise rather than on average), we recover the notion of a relatively measure-preserving extension.

Our motivation is threefold. First, we develop a theory of relative stationary dynamical systems, which generalises the well-studied concept of stationary dynamical systems. Within this theory, we prove an analogous structure theorem in the sense of Furstenberg-Glasner~\cite{GF09}. Furstenberg-Glasner \cite{GF09} asked whether stationary dynamical systems can always be naturally related to a factor of some universal object. Their Structure Theorem \cite[Theorem 4.3]{GF09} provides an affirmative answer to this question.  
We prove that a stationary extension can be similarly related to a factor of some space which is universal in an appropriate sense. 
This leads us to the notion of \textit{Relative Poisson Boundary}, our second motivation. Just as the Poisson boundary is the probabilistic boundary for the group, there should exist a unique maximal extension of a base space $(Y,\eta)$ that captures the relative asymptotic behaviour of random walks.
Third, this work is motivated by recent developments in the stationary theory of operator algebras, specifically the work of Hartman and Kalantar \cite{HartKal}. Relative stationarity can be expressed in terms of conditional expectations. This opens new paths to be explored in the non-commutative setup.

Our first main result establishes the existence and uniqueness of the maximal relative boundary.

\begin{theorem}[Relative Poisson Boundary]
\label{thm:Poissonrel_intro}
Let $(Y,\eta)$ be a standard non-singular $G$-space and let $\mu$ be an admissible non-degenerate Borel probability measure on $\G$. Let $\poi(\G,\mu)$ denote the standard Poisson boundary of $(\G,\mu)$. Then the projection map 
\[ \pr_2: \poi(\G,\mu)\otimes (Y,\eta)\longrightarrow (Y,\eta) \] 
is the maximal relative $\mu$-boundary over $(Y,\eta)$. Furthermore, it is unique up to $G$-isomorphism.
\end{theorem}

This theorem provides a concrete realisation of the relative boundary; it is the product of the Poisson boundary with the base space, though the measure on the product need not be a product measure in the general case (though it is for the maximal boundary). This generalises the intuition that relative boundaries should capture \say{fiber-wise} randomness.

Our second main contribution is the Relative Structure Theorem. To do this, we use a notion of relative joining over a measurable space (see Section~\ref{sec: join} for a notion of relative joining $\nu\vee \xi$) and introduce the notion of being relatively $1-1$ in the limit (see Definition~\ref{def: rel 1-1}).

\begin{theorem}[Relative Structure Theorem]\label{thm: structure intro}
Let $\G$ be a locally compact second countable group with an admissible non-degenerate probability measure $\mu$. Let $\pi: (X,\nu)\longrightarrow (Y,\eta)$ be a relatively $\mu$-stationary map. Then there exists a relative $\mu$-boundary map $\phi: (\prob(X),\xi)\longrightarrow (Y,\eta)$ such that the following diagram commutes:
$$      \xymatrix{ &(X\times \prob(X),\nu \vee \xi)\ar[dd]^{\psi}\ar[dl]^{\pr_1}\ar[dr]^{\pr_2} & \\
      (X,\nu)\ar[dr]^{\pi} & &(\prob(X),\xi)\ar[dl]^{\phi}
        \\ 
       & (Y,\eta)&
        }        $$
where the joining $\nu \vee \xi$ is constructed via boundary limits, $\pr_1$ is relatively 1-1 in the limit, and $\pr_2$ is relatively measure-preserving.
\end{theorem}

\subsection*{Organisation of the Paper} Apart from this section, this paper contains six other sections. We introduce the notations used throughout in Section~\ref{sec:notationandpreliminaries}. Moreover, we also prove some well-known facts about compact models and the disintegration map in this section. We introduce the notion of relative stationarity in Section~\ref{sec:relatstation} and prove some of its basic properties. In particular, we show that relative stationarity passes to intermediate factors. In Subsection~\ref{subsec:ucp}, we introduce the notion of stationarity for general unital completely positive maps in the context of $C^*$-dynamical systems. In Section~\ref{sec:relprox}, we introduce the notion of relative proximality. In particular, we show the existence of a measurable and $G$-equivariant relative boundary map (see Corollary~\ref{corollary: bnd to prob}). We also show that the identity map is the only endomorphism between relative $\mu$-boundaries (see Proposition~\ref{proposition: id map}). We introduce the notion of relative boundaries in Section~\ref{sec:relboundary}. Apart from giving various examples of relative boundaries, we also show that relative boundary is a hereditary property. In Section~\ref{sec: join}, we talk about relative joining over a measurable space $(Y,\eta)$ and prove the uniqueness of such a joining in the case of a relative boundary (see Proposition~\ref{proposition:uniqejoining}). We prove the Structure theorem in Section~\ref{sec: str}. Finally, in Section~\ref{Sec: Poissonbdryrelative}, we establish the existence and uniqueness of the relative Poisson boundary.
\subsection*{Acknowledgment}The first and last-named authors thank Mehrdad Kalantar and Yair Hartman for motivating the problem in the first place. Part of the work was done when the first-named author visited the Department of Mathematics at Universität Innsbruck. The author thanks the Department for their hospitality. We also thank David Beck-Tiefenbach for many valuable discussions. We are grateful to Yongle Jiang for taking the time to read through a near-complete draft of this paper and for his numerous suggestions and corrections. 
\subsection*{Funding} This research was supported by a grant from the Israel Science Foundation:
ISF 1175/18 for the first author. The first-named author is also supported by research funding from the European Research Council
(ERC) under the European Union's Seventh Framework Program
(FP7-2007-2013) (Grant agreement No. 101078193). This research was funded in whole or in part by the Austrian Science Fund (FWF) 10.55776/PAT3123425. The two last-named authors were partially supported by Early Stage Funding- Vice-rector for research, University of Innsbruck. 
The last-named author was also partly supported by the French national grant ANR-22-CE40-0004 of the project GoFR. 
This research was funded in whole or in part by the \textbf{Austrian Science Fund (FWF)
10.55776/ESP4189024}.
This work was partially supported by the Simons Foundation grant (award no. SFI-MPS-T-Institutes-00010825) and by State Treasury funds as part of a task commissioned by the Minister of Science and Higher Education under the project “Organization of the Simons Semesters at the Banach Center - New Energies in 2026-2028” (agreement no. MNiSW/2025/DAP/491).

\section{Notation and Preliminaries}
\label{sec:notationandpreliminaries}
Let $\G$ be a locally compact, second countable group with a probability measure $\mu$. We assume $\mu$ is admissible and non-degenerate, i.e., there exists a convolution power $k$ such that $\mu^{*k}$ is absolutely continuous with respect to the Haar measure $m$ on $\G$ and that the support of $\mu$ generates $\G$ as a semigroup. Throughout let  $(X,\nu)$, $(Y,\eta)$, $(Z,\zeta)$, $(W,\xi)$
denote standard Borel probability spaces on which $\G$ acts measurably and non-singularly (i.e. $\nu$ and $g\nu$ have the same null-sets, for all $g\in G$).
\begin{definition}[Compact Model]
Given $ \G\curvearrowright(X,\nu)$, a compact model for $\G\curvearrowright(X,\nu)$ is a compact Hausdorff space $X_0$ such that the action $G\curvearrowright X_0$ is continuous and  $C(X_0)\subset L^{\infty}(X,\nu)$ is dense in the weak$^*$-topology. 
\end{definition}
The existence of a compact model is well-known and can be traced back to \cite{Var63} (also see \cite[Proposition~2.1]{amrutam2024subalgebras}). Nonetheless, we shall show its existence (in a much more general setup) below. Consequently, throughout this note, we will pass to compact models and, without any loss of generality, assume that the spaces $X, Y, Z$ we deal with are compact and the action is continuous. We call them $\G$-spaces.

A \textit{$\G$-factor map} or just \textit{factor map} $\pi:(X,\nu)\longrightarrow (Y,\eta)$ is a measurable surjective map $\pi: X_0\longrightarrow Y_0$ such that $\pi\nu=\eta$ and $\g\pi(x)=\pi(\g x)$, for all $\g\in\G$, $x\in X_0$ (\textit{$\G$-equivariance}), where $X_0$ and $Y_0$ are measurable sets with $ X_0=X$ mod $\nu$ and $Y_0=Y$ mod $\eta$. 

\begin{proposition}
\label{compactmodelsforfactormaps}
Let $\pi:(X,\nu)\longrightarrow (Y,\eta)$ be a $G$-factor map between non-singular standard Borel probability spaces. We can find compact models $X_0$ (for $X$) and $Y_0$ (for $Y$) and a map $\rho: X_0\to Y_0$ such that $\rho$ is continuous. Moreover, if we let $\pi_{X_0}$ and $\pi_{Y_0}$ be the corresponding maps from $X\to X_0$ and $Y\to Y_0$ respectively, then $\pi_{Y_0}^{-1}\circ\rho\circ\pi_{X_0}=\pi$.
\begin{proof}
We follow the arguments contained in \cite[Lemma~2.6]{creutz2014normal} and adapt them to our purposes (also see \cite[Appendix]{zimmerBook}). Since the unit ball of $L^{\infty}(X,\nu)$ is compact and metrizable in the weak$^*$-topology, we can find a countable dense collection of functions $\{f_n:n\in\mathbb{N}\}\subset L^{\infty}(X,\nu)$. Similarly, let $\{h_n: n\in\mathbb{N}\}$ be a countable dense collection of functions in $L^{\infty}(Y,\eta)$. Let $\mathcal{F}=\cup_{n\in\mathbb{N}}\left\{\{f_n\}\cup \{h_n\circ\pi\}\right\}$ and $\mathcal{F}'=\cup_{n\in\mathbb{N}}\{h_n\}$. We denote by $\mathbb{B}$ the unit ball of $L^{\infty}(G,m)$. Since $\mathbb{B}$ is compact in the weak$^*$-topology, the spaces $\prod_{f\in\mathcal{F}} \mathbb{B}$ and $\prod_{f\in\mathcal{F}'} \mathbb{B}$ are compact in the product topology. We consider $\prod_{f\in\mathcal{F}} \mathbb{B}$ as a $G$-space with respect to the right action, i.e., for $g\in G$ and $\left(\phi_f\right)_{f\in\mathcal{F}}\in \prod_{f\in\mathcal{F}} \mathbb{B}$,
\[g.\phi_f(h)=\phi_f(hg), \ h\in G,\  f\in\mathcal{F}.\]
Similarly, $\prod_{f\in\mathcal{F}'} \mathbb{B}$ is a $G$-space with respect to the right action.
Define $\Phi_X: X\to \prod_{f\in\mathcal{F}} \mathbb{B}$ by $\Phi_X(x)=\left(\phi_f(x)\right)_{f\in\mathcal{F}}$, where $\phi_f(x)(g)=f(gx)$. Similarly, define $\Phi_Y: Y\to \prod_{h\in\mathcal{F}'} \mathbb{B}$ by $\Phi_Y(y)=\left(\phi_h(y)\right)_{h\in\mathcal{F}'}$, where $\phi_h(y)(g)=h(gy)$. Now, for $h\in G$, $f\in\mathcal{F}$,
\[\phi_f(hx)(g)=f(ghx)=(\phi_f(x))(gh)=h.(\phi_f(x))(g),\ g\in G.\]
This shows that $\Phi_X$ and $\Phi_Y$ are $G$-equivariant. Moreover, let $\pi_0:\prod_{f\in\mathcal{F}} \mathbb{B}\to \prod_{f\in\mathcal{F}'} \mathbb{B}$ be the restriction map, i.e.,
for $f\in\mathcal{F}'$, define the $f$\textsuperscript{th} coordinate of $\pi_0(x)$ as $(f\circ\pi)$\textsuperscript{th} coordinate of $x$. Since $\mathcal{F}$ and $\mathcal{F}'$ form a separating family for $L^{\infty}(X,\nu)$ and $L^{\infty}(Y,\eta)$ respectively, the maps $\Phi_X$ and $\Phi_Y$ are injective maps. Let $X_0=\overline{\Phi_X(X)}$ and $Y_0=\overline{\Phi_Y(Y)}$. We now have the following diagram of commutative maps:
\begin{center}
    \begin{tikzcd}
(X,\nu) \arrow{d}{\Phi_X} \arrow{r}{\pi}
& (Y,\eta) \arrow{d}{\Phi_Y}\\
(X_0,(\Phi_X)_*\nu) \arrow{r}{\pi_0}
& (Y_0,(\Phi_Y)_*\eta)
\end{tikzcd}.
\end{center}
Since $\Phi_X$ and $\Phi_Y$ are $G$-isomorphisms onto $X_0$ and $Y_0$, letting $\Pi_{X_0}=\Phi_X$, $\rho=\pi_0$ and $\Pi_{Y_0}=\Phi_Y$, the claim follows.
\end{proof}
\end{proposition}  

\begin{definition}[Disintegration of measures]
Let $\pi:(X,\nu)\longrightarrow(Y,\eta)$ be a factor map. Then 
for $\eta$-almost every $y\in Y$ there exists a probability measure $\nu_y$ on the fibre $\pi^{-1}(y)$ such that
$$ y\mapsto \nu_y(A\cap \pi^{-1}(y)) \text{ is measurable}$$ for all measurable subsets $A\subseteq X$,
and 
$$ \int_{Y} \nu_y\, d\eta(y)=\nu.$$
The collection of probability measures $\{\nu_y\}_{y\in Y_0}$ is unique with these properties and is called \textit{disintegration of $\nu$ w.r.t. $\pi$}. The existence and uniqueness were proven by Rokhlin in \cite{Rokhlin}. When $X$ is a compact $G$-space, the condition $ y\mapsto \nu_y(A\cap \pi^{-1}(y))$  being measurable is equivalent to saying that the map $ y\mapsto \nu_y(f)$ is measurable for every $f\in C(X)$.
\end{definition}

%The factor map $\pi$ is called \textit{relatively measure-preserving} if $\nu_{\g y}=\g\nu_{y}$ for all $\g\in\G$ and $\eta$-a.e. $y\in Y$. 

Let $X$ be a compact metric space on which $\G$ acts continuously and let $\mathcal{B}(X)$ denote the Borel $\sigma$-algebra on $X$. Then, there is a natural continuous action of $\G$ on  $\prob(X)$, the space of Borel probability measures on $X$ equipped with the weak*-topology. This action is defined by the push-forward measure
$$\g\nu(A):=\nu(g^{-1}(A)), \ \forall A\in \mathcal{B}(X) $$
for $\g\in \G$ and  $\nu\in\prob(X)$.
Recall that the weak*-topology is the initial topology on $\text{Prob}(X)$ induced by the continuous functions $f\in C(X)$, i.e., the coarsest topology on $\text{Prob}(X)$ such that for each $f\in C(X)$, the map $\Phi_f:\text{Prob}(X)\to\mathbb{C}$, $\nu\mapsto\nu(f)$ is continuous. 

It is well-known to the experts that the disintegration map is measurable. We include the proof nonetheless for the sake of completeness. 
\begin{proposition}
\label{proposition:measurabilityofthedisint}
Let $X$ be a compact Hausdorff space on which $G$ acts continuously.
Let $\pi:(X,\nu)\longrightarrow(Y,\eta)$ be a factor map. Let $D_{\pi}: Y\to\text{Prob}(X)$, $y\mapsto\nu_y$ be the $\eta$-almost everywhere defined disintegration map. Then, $D_{\pi}$ is measurable w.r.t. the weak*-topology on $\prob(X)$.
\begin{proof}
Since the topology on $\text{Prob}(X)$ is the initial topology induced by the continuous functions $f\in C(X)$, the topology is generated by open sets of the form 
\[\left\{\Phi_f^{-1}(U): f\in C(X), U\subseteq\mathbb{C}\text{ open}\right\}.\]
It is enough to show that $D_{\pi}^{-1}\left(\Phi_f^{-1}(U)\right)$ is measurable for each $f\in C(X)$ and $U\subseteq\mathbb{C}$ open. Since $$D_{\pi}^{-1}\left(\Phi_f^{-1}(U)\right)=\left(\Phi_f\circ D_{\pi}\right)^{-1}(U)=\{y\in Y: \nu_y(f)\in U\},$$
and the latter set $\{y\in Y: \nu_y(f)\in U\}$ is measurable by definition of the disintegration, the claim follows.
\end{proof}
\end{proposition}
%$G\curvearrowright C(Y)$ on $C(Y)$ defined by $gf(x)=f(g^{-1}x)$. Moreover, $G\curvearrowright P(Y)$, the set of probability measures on $Y$ defined by \[g\eta(f)=\eta(g^{-1}f)=\int_Yg^{-1}fd\eta=\int_Y f d(g\eta).\]

\begin{definition}\label{def: rel mp}
    A $\G$-factor map $\pi: (X,\nu)\longrightarrow (Y,\eta)$ is called \textit{relatively measure-preserving} if 
    $$\g\nu_{\g^{-1}y}=\nu_y,
    $$ for all $\g\in \G$ and $\eta$-almost all $y\in Y$.
\end{definition}

\section{Relative Stationarity}
\label{sec:relatstation}
We now introduce a generalisation of the concept of measure-preserving maps. % from Definition~\ref{def: rel mp}. 
Namely, we call a map relatively stationary if the equation $\g\nu_{\g^{-1}y}=\nu_y    $ holds ``on average'' with respect to a given probability measure $\mu$ on the group $\G$. This is in the same spirit as generalising measure-preserving systems to $\mu$-stationary systems. The only difference is that here, we consider relations between systems rather than a single system.

\begin{definition}[Relatively $\mu$-stationary]
\label{def:stationary}
We say that a $\G$-factor map  $\pi:(X,\nu)\longrightarrow (Y,\eta)$ is \textit{relatively $\mu$-stationary} if 
$$\nu_y=\int_{\G} g\nu_{\g^{-1}y}\, d\mu(\g), \ \text{ for $\eta$-a.e. $y\in Y$.}
$$
\end{definition}
Note that in the case of  $(Y,\eta)$ being a one-point space with a trivial $\G$-action, the above definition is equivalent to $\G\curvearrowright(X,\nu)$ being $\mu$-stationary in the usual sense. This is proven in the following proposition, together with other basic properties of our notion.

% Let  $\g\in \G$ be given. Consider $\widetilde{\pi}:(X,\g \nu)\longrightarrow (Y,\g \eta)$. Let $\{\nu_y\}$ be the disintegration w.r.t. $\pi$. Then     $\{\g\nu_{\g^{-1} y}\}$ is the disintegration w.r.t. $\widetilde{\pi}$. Indeed, we see that     $$\int_Y g\nu_{\g^{-1} y}\, d\g\eta(y)=\int_Y \g \nu_{ y}\, d\g\eta(\g y)    =\int_Y \g\nu_{ y}\, d\eta( y)=\g \nu.    $$

\begin{proposition}[Basic properties]\label{proposition: basic} Let  $\pi:(X,\nu)\longrightarrow (Y,\eta)$ be a $\G$-factor map and $\mu$ be a Borel probability measure on $\G$. 
\begin{enumerate}
    \item\label{basic: rel mp} Every relatively measure-preserving factor map is relatively $\mu$-stationary.
    \item\label{basic: inv}  Assume that $\eta$ (or $\nu$) is invariant under elements of the support of $\mu$. Then $\pi$ is relatively $\mu$-stationary if and only if the measure $\nu$ is $\mu$-stationary.\\
    In particular, $
    \pi:(X,\nu)\longrightarrow (\{\star\}, \delta_{\star})$  is relatively $\mu$-stationary if and only if $\nu$ is $\mu$-stationary.
     %\item\label{basic: proj} The projection $\pr_2:(X,\nu)\otimes (Y,\eta)\longrightarrow (Y,\eta)$ is relatively $\mu$-stationary if and only if $\nu$ is $\mu$-stationary in the usual sense.
     \item\label{basic: stat eta} If $\eta$ is $\mu$-stationary and $\pi$ is relatively measure preserving, then $\nu$ is $\mu$-stationary.
\end{enumerate}
  \end{proposition}

\begin{proof}
    The first claim is obvious since $\g\nu_{\g^{-1}y}=\g\g^{-1}\nu_{y}=\nu_y$ in case $\pi$ is relatively measure-preserving. \\
     To show (\ref{basic: inv}) let us assume that $\g \eta=\eta$ for all $\g\in\supp(\mu)$. In this case, we obtain
    $$\int_{Y} \int_{\G} \g\nu_{\g^{-1} y}\,  d\mu(g)\,  d\eta(y)=    
    \int_{Y} \int_{\G} \g\nu_{ y}\,  d\mu(g)\, d\eta(y)= 
    \int_{\G} \g\nu\,  d\mu(g)  $$ 
    and thus $\int_{\G} \g\nu_{\g^{-1} y}\,  d\mu(g)=\nu_y$  if and only if $
    \int_{\G} \g\nu\,  d\mu(g)=\nu$, using the uniqueness of the disintegration for the if-direction.\\
    We now show (\ref{basic: stat eta}).
    \begin{align*}
    \int_{\G} \g\nu\,  d\mu(g)&=
    \int_{\G} \int_{Y}(\g \nu)_{y}\, d\eta(y)\, d\mu(\g)\\&=
       \int_{\G} \int_{Y} \g \nu_{ y}\, d\eta(y)\, d\mu(\g)=  \int_{\G} \int_{Y} \nu_{\g y}\, d\eta(y)\, d\mu(\g)
       \\&=
       \int_{\G} \int_{Y} \nu_{ y}\, d\g\eta(y)\, d\mu(\g)
      \\& =       \int_{Y} \nu_{ y}\, d\eta(y) =\nu\end{align*}
     due to relatively measure-preservingness and $\mu$-stationarity of $\eta$. 
\end{proof}

It is a well-known fact that if $G\curvearrowright(X,\nu)$ is ergodic and $\mu$-stationary and $G$ is countable, then either $\nu$ is atomless (i.e. $\nu(\{x\})=0,$ $\forall x\in X$), or $\nu$ is $G$-invariant. % and $X$ is finite. %-%$X$ finite is true, but has nothing to do with stationarity%-%
We get a similar result in the relative $\mu$-stationary setup as follows.

\begin{proposition}\label{prop: atomless or rel mp}
    Let $G$ be countable and let $G \curvearrowright(X,\nu)$ be ergodic and non-singular. If $\pi:(X,\nu)\longrightarrow (Y,\eta)$ is relatively $\mu$-stationary, then either $(X,\nu)$ is atomless, or $\pi$ is relatively measure preserving. % and $X$ is finite.
\end{proposition}
\begin{proof} W.l.o.g. we may assume that $X$ is compact.
Assume there is $w\in X$ such that $\nu(\{w\})>0$. Then, by ergodicity, $Gw=X$ mod $\nu$. Thus, $X$ is countable up to a null-set $N$. Moreover, since the action of $G$ on $X\setminus N$ is transitive, we get that $X\setminus N$ is discrete. Hence $N$ is open and so $X\setminus N$ is compact again. We conclude that $X\setminus N$ has to be finite. %Since $X$ is assumed to be compact, we obtain that $X$ is finite up to a null set. 
%Thus, we can pick an element $z\in X$ such that $\nu(\{z\})>0$ is maximal among all atom-values $\nu(\{x\})$ for $x\in X$. (There might not be a unique such $z$, but we can always find at least one.) 
Note that $\eta(\{\pi(x)\})=\nu(\pi^{-1}(\{\pi(x)\})) \geq \nu(\{x\})$, therefore $\nu_{\pi(x)}$ is defined, whenever $\nu(\{x\})>0.$ 
Since $x$ lies in exactly one fieber, namely $\pi^{-1}(\{\pi(x)\})$, %for $\nu(\{x\})>0$,
whenever $\nu_{\pi(x)}$ is defined, we get that
$$\nu(\{x\})=\int_Y \nu_y(\{x\})\, d\eta(y)=\nu_{\pi(x)}(\{x\})\eta(\{\pi(x)\}),$$
because $\supp(\nu_y)\subseteq \pi^{-1}(\{y\})$. One can see that  $\nu(\{x\})>0$ if and only if $\nu_{\pi(x)}$ is defined and  $\eta(\{\pi(x)\})$ and $\nu_{\pi(x)}(\{x\})$ are both strictly positive. 
% -% Since $X$ is finite up to a null set, also $Y$ is finite up to a null set. Indeed, $\pi^{-1}(Y)$ is finite mod $\nu$, and a factor map is understood in the measurable category, hence we identify factor maps which only differ on a null-set. Hence, we can think of $ \pi$ as being only defined on the finite set. Since $\pi$ is surjective, we obtain that $Y$ is finite. %-%
%Let $N$ be a measurable $\nu$-null-set such that $X\setminus N$ is finite.
Thus, the set $$S:=\{\nu_{\pi(x)}(\{x\}) \, : \, x\in X\setminus N, \, \eta(\{\pi(x)\})>0\}$$ is finite and non-empty, and thus has a maximal value, which is strictly positive. Choose one $z\in X \setminus N$ such that $\nu_{\pi(z)}(\{z\})$ is maximal. (There might not be a unique such $z$, but we can always find one.) 
Now, by relative $\mu$-stationarity, we have
$$\nu_{\pi(z)}(\{z\})=\sum_{g\in G}\mu(g) g\nu_{g^{-1}\pi(z)}(\{z\})
=\sum_{g\in G}\mu(g) \nu_{\pi(g^{-1}z)}(\{g^{-1}z\})
$$
But the left-hand side is a convex combination of elements in $S$. Indeed, $\eta(\{ \pi(g^{-1}z)\})=\eta(g^{-1}\{\pi(z)\})>0 $ by non-singularity and the assumption that $\eta(\pi(\{z\}))>0$. And clearly $g^{-1}z\in X\setminus N=Gw$, since $z\in Gw$.
 By maximality of $\nu_{\pi(z)}(\{z\})$, we thus obtain that 
 $$
     g\nu_{g^{-1}\pi(z)}(\{z\})= \nu_{\pi(z)}(\{z\}), \ \forall g\in\supp(\mu).
      $$
 The above extends to all $g\in G$, since $\mu$ is non-degenerate.
Moreover, since $X=Gz=Gw$ mod $\nu$, we know that for a.e. $x\in X$ there exists $g\in G$ such that $x=gz$. Therefore, the above equation holds for almost all points in $X$, not only the particularly chosen $z$, i.e.  
\begin{equation}\label{eq: rel mp on singelton}
 \nu_{\pi(x)}(\{x\})= \nu_{g\pi(x)}(g\{x\}) \ \text{ for $\nu$-a.e. $x\in X$}, \ \forall g\in G.
\end{equation} 
Fix $x\in X\setminus N$. It is left to show that the above equality still holds when replacing $\{x\}$ by any measurable set in the fieber $\pi^{-1}(\{x\})$.
 Let $A\subseteq \pi^{-1}(\{\pi(x)\})$ be measurable. 
 Then there exists $F_A\subseteq G$ finite, such that $A\cap N^c=\bigcup_{g\in F_A} \{gx\}$ because $Gx = Gw=X\setminus N$. W.l.o.g. we may assume that $gx\neq x$ for all $g\in F_A\setminus\{e\}$.
 Thus, whenever $\nu_{\pi(x)}$ is defined, we can write
 $$
 \nu_{\pi(x)}(A\cap N^c)=\sum_{g\in F_A} \nu_{\pi(x)}(\{gx\}).
 $$
Now, since $gx\in A\subseteq \pi^{-1}(\{\pi(x)\})$, we see that $\pi(gx)=\pi(x)$ for all $g\in F_A$. Thus, the above  expression becomes
$$
\nu_{\pi(x)}(A\cap N^c)=\sum_{g\in F_A} \nu_{\pi(gx)}(\{gx\})=\nu_{\pi(x)}(\{x\}) \cdot\vert F_A \vert,
$$
by Equation~\ref{eq: rel mp on singelton}.
Therefore, for every given $g\in G$, we obtain
$$
\nu_{\pi(x)}(A\cap N^c)= \nu_{\pi(x)}(\{x\}) \cdot\vert F_A \vert= \nu_{\pi(gx)}(\{gx\}) \cdot\vert F_{gA} \vert= \nu_{g\pi(x)}(gA\cap N^c),
$$
again using Equation~\ref{eq: rel mp on singelton} and the obvious fact that $\vert A\vert=\vert gA\vert$.
Since $N$ is a null-set, the above equation proves that $\pi$ is relatively measure preserving.
 \iffalse %------------------------------------
Note that $\nu_{\pi(x)}(\{gz\})=0$ whenever $\pi(x)\neq \pi(gz)$. By construction, there can be at most one $g\in F_A$ such that $g z=x$, hence $$
    \nu_{\pi(x)}(A\cap N^c)= 
    \begin{cases} \nu_{\pi(x)}(\{x\}), &\text{ if } x \in F_A z,\\
    0, & \text{ else }.
\end{cases}
$$
Now, by Equation~\ref{eq: rel mp on singelton}, 
we obtain for a.e. $x\in X$, and whenever $\nu_{\pi(x)}$ is defined, that
\begin{align*}
 \nu_{\pi(x)}(A\cap N^c)=&
  \nu_{\pi(x)}(\{x\}) \chi_{x\in F_A z} 
 \\
 =&  \nu_{\widetilde{g}\pi(x)}(\widetilde{g}\{x\})  \chi_{x\in F_A z}
 \\
 =&  \nu_{\widetilde{g}\pi(x)}(\widetilde{g} (A\cap N^c))
\end{align*}
for all $\widetilde{g}\in G$.
Since $\nu(\cdot )=\nu(\cdot \cap N^c)$, by uniqueness of the disintegration, we have $\nu_y(\cdot \cap N^c)=\nu_y(\cdot)$ for $\eta$-a.e. $y\in Y$.
We can thus conclude from the above that
$$ \nu_{\pi(x)}(A)= \nu_{\pi(x)}(A\cap N^c)=  \widetilde{g}^{-1}\nu_{\widetilde{g}\pi(x)}(A\cap N^c)=  \widetilde{g}^{-1}\nu_{\widetilde{g}\pi(x)}( A),
$$ whenever the expression is defined.
This shows that $$ \nu_{\pi(x)}= \widetilde{g}^{-1}\nu_{\widetilde{g}\pi(x)}, \ \forall \widetilde{g}\in G, \text{ for a.e. } x\in X,
$$ i.e. $\pi$ is relatively measure-preserving.
\fi %-----------------------------
\end{proof}

\begin{remark}
   We believe that there are factor maps $\pi:(X,\nu)\longrightarrow (Y,\eta)$ that are not relatively $\mu$-stationary, even when $\nu$ is $\mu$-stationary in the usual sense.  In the case where $(Y,\eta)$ is ergodic and $G$ is countable,  the conjectured phenomenon can only happen in the case when $Y$ is uncountable. Because, otherwise $\eta$ would be invariant and thus Proposition~\ref{proposition: basic} (\ref{basic: inv}) would give that $\pi$ is relatively $\mu$-stationary.
   By Proposition~\ref{prop: atomless or rel mp} is sufficies to find a map which is not relatively measure preserving but there are atoms in $(X,\nu)$.
\end{remark}
\begin{example}
Let $\pi: (X,\nu)\to (Y,\eta)$ be a $G$-factor map. Assume that $G=\mathbb{Z}/2\mathbb{Z}=\{e,s\}$. Let $Y=\{y_1,y_2\}$ with $\eta(y_i)\ne \{0\}$. We also assume that $sy_1=y_2$ and consequently, $sy_2=y_1$. Also assume that $\nu$ is $\mu$-stationary for some $\mu\in\text{Prob}(G)$. We show that $\pi$ is relatively $\mu$-stationary. 
Combining the $\mu$-stationarity of $\eta$ along the fact that $\eta\in\text{Prob}(Y)$, we see that 
\[\mu(e)(1-\eta(y_2))+\mu(s)\eta(y_2)=1-\eta(y_2).\]
This implies that 
\[\eta(y_2)=\frac{1-\mu(e)}{1+\mu(s)-\mu(e)}=\frac{\mu(s)}{2\mu(s)}=\frac{1}{2}.\]
Therefore, $\eta(y_1)=\frac{1}{2}$. This shows that $\eta$ is an invariant measure on $Y$.
Since $\nu_{y_1}\eta(y_1)+\nu_{y_2}\eta(y_2)=\nu$, we see that
\begin{align*}
&\mu(e)\nu_{y_1}\eta(y_1)+\mu(e)\nu_{y_2}\eta(y_2)=\mu(e)\nu, \text{and}\\&\mu(s)\nu_{y_1}\circ s^{-1}\eta(y_1)+\mu(s)\nu_{y_2}\circ s^{-1}\eta(y_2)=\mu(s)\nu\circ s^{-1}
\end{align*}
Since $sy_1=y_2$, adding these two equations together along with the fact that $\nu$ is $\mu$-stationary, we see that
\[\left(\mu(e)\nu_{y_1}\eta(y_1)+\mu(s)\nu_{sy_1}\circ s^{-1}\eta(y_2)\right)+\left(\mu(e)\nu_{y_2}\eta(y_2)+\mu(s)\nu_{sy_2}\circ s^{-1}\eta(y_1)\right)=\nu\]
Since $\eta(y_1)=\eta(y_2)$, this further implies that
\[\left(\mu(e)\nu_{y_1}+\mu(s)\nu_{sy_1}\circ s^{-1}\right)\eta(y_1)+\left(\mu(e)\nu_{y_2}+\mu(s)\nu_{sy_2}\circ s^{-1}\right)\eta(y_2)=\nu\]
Since the probability measures $\mu(e)\nu_{y_1}+\mu(s)\nu_{sy_1}\circ s^{-1}$ and $\mu(e)\nu_{y_2}+\mu(s)\nu_{sy_2}\circ s^{-1}$ are supported on $\pi^{-1}(\{y_1\})$ and $\pi^{-1}(\{y_2\})$ respectively, by the uniqueness of the disintegration measures, it follows that
\[\mu(e)\nu_{y_1}+\mu(s)\nu_{sy_1}\circ s^{-1}=\nu_{y_1}\text{ and }\mu(e)\nu_{y_2}+\mu(s)\nu_{sy_2}\circ s^{-1}=\nu_{y_2}.\]
\end{example}
\subsection*{Stationarity of projection maps}
Given two $G$-spaces $(X,\nu)$ and $(Y,\eta)$, the product space $(X,\nu)\otimes (Y,\eta)$ along with the projection constitute a class of examples which have been studied a lot in the context of intermediate factor theorems (see for example~\cite{BS06}).   
We consider this setup below. 
  \begin{example}\label{ex: proj} The projection $\pr_2:(X,\nu)\otimes (Y,\eta)\longrightarrow (Y,\eta)$ is relatively $\mu$-stationary if and only if $\nu$ is $\mu$-stationary in the usual sense.
  Indeed, in this case, the disintegration measures are given by $(\nu\otimes \eta)_y=\nu\otimes \delta_{y}$ and thus $$\int_{\G} \g (\nu\otimes \eta)_{\g y}\, d\mu(\g)= \Big(\int_{\G} \g \nu \, d\mu(\g) \Big)\otimes\delta_{y}=\nu\otimes \delta_{y}$$ where the last equality holds if and only if $\nu$ is $\mu$-stationary.
\end{example}
\subsection*{Heredity to intermediate factors}
The pushforward of a stationary measure is stationary. Similarly, the relative $\mu$ stationarity behaves well with respect to factors.
\begin{proposition}(Intermediate factors)\label{proposition: IF}
        Consider a relatively $\mu$-stationary factor map $\pi: (X,\nu)\longrightarrow (Y,\eta)$. 
        Let $(Z,\eta)$ be a $\G$-space such that $$         \xymatrix{ (X,\nu)\ar[dd]^{\pi}\ar[dr]^{\kappa} & 
        \\ 
        & (Z,\zeta)\ar[dl]^{\sigma}\\
        (Y,\eta)
        }        $$ 
        is a
        commutative diagram of factor maps. Then $\sigma: (Z,\zeta)\longrightarrow (Y,\eta)$ is relatively $\mu$-stationary. 
\end{proposition}

\begin{proof}
First note that $\zeta_{y}=\nu_{y}\circ \kappa^{-1}$ for $\eta$-almost every $y\in Y$, meaning that $\zeta_z(\sigma^{-1}(y)\cap A)=\nu_y(\pi^{-1}(y)\cap \kappa^{-1}(A))$ for all $A\in\mathcal{B}(Z)$. Indeed, $\int \nu_y(\pi^{-1}(y)\cap \kappa^{-1}(A))\, d\eta(y)=\nu(\kappa^{-1}(A))=\zeta(A)$ and thus the claim follows by the uniqueness of the disintegration. Now, for $\eta$-almost every $y\in Y$ we obtain $$\int_{\G} \g\zeta_{\g^{-1}y}(A\cap \sigma^{-1}(\g^{-1}y))\, d\mu(\g)
=
\int_{\G} \g\nu_{\g^{-1}y}(\kappa^{-1}(A)\cap \pi^{-1}(\g^{-1}y))\, d\mu(\g)
$$ 
$$=\nu_{y}(\kappa^{-1}(A))=\zeta_y(A)
$$ by relative $\mu$-stationarity of $\pi$. This proves that $\sigma$ is relatively $\mu$-stationary.    
\end{proof}
In general, it is not enough that $\kappa$ and $\sigma$ are relatively $\mu$-stationary to conclude that $\pi$ is relatively $\mu$-stationary in the above picture, as the following easy example demonstrates.
\begin{example}
Let $(X,\nu)$ and $(Y,\eta)$ be  $(G,\mu)$-spaces. Consider the following diagram of $\G$-factor maps:
$$         \xymatrix{ (X\times Y,\nu\otimes\eta)\ar[dd]^{\pi\circ\pr_2}\ar[dr]^{\pr_2} & 
        \\ 
        & (Y,\eta)\ar[dl]^{\pi}\\
        \{\star\}
        }        .$$ 
Since $\eta$ is a $\mu$-stationary measure, $\pi$ is a relatively $\mu$-stationary map. It follows from Example~\ref{ex: proj} that the projection $\pr_2$ is a relatively $\mu$-stationary factor map. However, $\nu\otimes\eta$ may not be a $\mu$-stationary measure.     
\end{example}

\subsection{Existences and examples for ucp maps}
\label{subsec:ucp}
Even though the notion of relative stationarity is defined for the disintegration map, the definition can be adapted to a much more general framework. Just to briefly recall, a unital $C^*$-algebra $\mathcal{A}$ is a norm closed $*$-subalgebra of $\mathbb{B}(\mathcal{H})$ for some Hilbert space $\mathcal{H}$. A von Neumann algebra $\mathcal{M}$ is a weak operator topology (WOT)-closed  $*$-subalgebra of $\mathbb{B}(\mathcal{H})$.

Let us consider a discrete group $\Gamma$.  A unital $\Gamma$-$C^{\star}$-algebra $\mathcal{A}$ is a $C^{\star}$-algebra equipped with an action by $\Gamma$ through $\star$-automorphisms, ensuring that the operation $(g, a)\mapsto g.a$ is continuous for $g$ in $\Gamma$ and $a$ in $\mathcal{A}$. Similarly, a $\Gamma$-von Neumann algebra $\mathcal{M}$ is a von Neumann algebra with an action by $\star$-automorphisms. This action ensures that the operation mapping $(g, x)$ to $g.x$ for $g$ in $\Gamma$ and $x$ in $\mathcal{M}$ is continuous, similar to the case with the $\Gamma$-$C^{\star}$-algebra $\mathcal{A}$.
\begin{definition} 
\label{def:condstationary}
Let $\mu\in\text{Prob}(\Gamma)$. Given an inclusion of unital $\Gamma$-$C^*$-algebras $\mathcal{B}\subset\mathcal{A}$ with $\mathbf{1}_{\mathcal{B}}=\mathbf{1}_{\mathcal{A}}$, a unital completely positive map $\varphi:\mathcal{A}\to\mathcal{B}$ is called $\mu$-stationary if 
\begin{equation}
\label{eq:stationary}
\sum_{s\in\Gamma}\mu(s)s\varphi(s^{-1}a)=\varphi(a), \ a\in\mathcal{A}.\end{equation}
Similarly, given an inclusion of $\Gamma$-von Neumann algebras $\mathcal{N}\subset\mathcal{M}$, a normal unital completely positive map $\psi:\mathcal{M}\to\mathcal{N}$ is called $\mu$-stationary if satisfies equation~\eqref{eq:stationary} for every $a\in\mathcal{M}$.   
\end{definition}
If $\mathcal{M}=L^{\infty}(X,\nu)$ and $\mathcal{N}=L^{\infty}(Y,\eta)$ for some $G$-spaces $(X,\nu)$ and $(Y,\eta)$, then the canonical conditional expectation $\mathbb{E}_Y:L^{\infty}(X,\nu)\to L^{\infty}(Y,\eta)$ is defined by 
\[\mathbb{E}_Y(f)(y)=\nu_y(f), y\in Y_0,\]
where $Y_0\subset Y$ is a co-null set. In this context, to say that $\mathbb{E}_Y$ satisfies equation~\eqref{eq:stationary} is equivalent to saying that the disintegration measures are $\mu$-stationary in the sense of Definition~\ref{def:stationary}. In other words, an equivalent formulation of the factor map $\pi: (X,\nu)\to (Y,\eta)$ being relatively $\mu$-stationary is to say that the canonical conditional expectation $\mathbb{E}_Y: L^{\infty}(X,\nu)\to L^{\infty}(Y,\eta)$ is $\mu$-stationary in the sense of Definition~\ref{def:condstationary}.
\begin{example}
Let $\mathcal{M}$ be a $\Gamma$-von Neumann algebra represented faithfully inside $\mathbb{B}(\mathcal{H})$ for some Hilbert space $\mathcal{H}$.  
Denote by $\ell^2(\Gamma,\mathcal{H})$, the space of square summable $\mathcal{H}$-valued functions on $\Gamma$.
Note that $\Gamma\curvearrowright \ell^2(\Gamma,\mathcal{H})$ by left translation:
\[\lambda_g\xi(h):=\xi(g^{-1}h), \ \xi \in \ell^2(\Gamma,\mathcal{H}),\  g,h \in \Gamma.\]
Let $\sigma:\mathcal{M} \to B(\ell^2(\Gamma,\mathcal{H}))$ be the $*$-representation defined by \[\sigma(a)(\xi)(h):=\pi(h^{-1}a)\xi(h), \ a \in \mathcal{M}\]
where $\xi \in \ell^2(\Gamma,\mathcal{H})$, $h \in \Gamma$. The von Neumann crossed product
$\mathcal{M}\rtimes\Gamma$ is generated (as a von Neumann algebra inside $\mathbb{B}(\ell^
2
(\Gamma, \mathcal{H})$), by the left regular representation $\lambda$ of $\Gamma$ and the faithful $*$-representation $\sigma$ of $\mathcal{M}$ in $\mathbb{B}(\ell^2
(\Gamma, \mathcal{H}))$. There is a $\Gamma$-equivariant faithful normal conditional expectation $\mathbb{E}:\mathcal{M}\rtimes\Gamma\to\mathcal{M}$ defined by 
\[\mathbb{E}\left(\sigma(a_g)\lambda_g\right)=\left\{ \begin{array}{ll}
0 & \mbox{if $g\ne e$}\\
\sigma(a_e) & \mbox{otherwise}\end{array}\right\}.\]   Since $\mathbb{E}$ is $\Gamma$-equivariant, it follows that $\mathbb{E}$ is a $\mu$-stationary conditional expectation for any $\mu\in\text{Prob}(\Gamma)$. Of course, the conditional expectation associated with a $C^*$-algebraic crossed product $\mathcal{A}\rtimes_r\Gamma$ is also $\mu$-stationary for any $\mu\in\text{Prob}(\Gamma)$.
\end{example}
Given two von Neumann algebras $\mathcal{M}\subset\mathbb{B}(\mathcal{H})$ and $\mathcal{N}\subset\mathbb{B}(\mathcal{K})$, the von Neumann algebra $\mathcal{M}\overline{\otimes}\mathcal{N}$ is the von Neumann algebra generated by the algebraic tensor product $\mathcal{M}\otimes\mathcal{N}$ inside $\mathbb{B}(\mathcal{H}\otimes\mathcal{K})$.
\begin{example}
Consider the tensor product of two $\Gamma$-von Neumann algebras $\mathcal{M}\overline{\otimes}\mathcal{N}$. Let $\mu\in\text{Prob}(\Gamma)$ and $\varphi$ a $\mu$-stationary normal state on $\mathcal{M}$. Associate with $\varphi$, the right slice map $R_{\varphi}:\mathcal{M}\overline{\otimes}\mathcal{N}\to\mathcal{N}$, is defined by
\[R_{\varphi}(a\otimes b)=\varphi(a)b,~a\in\mathcal{M},~b\in\mathcal{N}.\]
We claim that $R_{\varphi}$ is $\mu$-stationary in the sense that 
\[\sum_{s\in\Gamma}\mu(s)sR_{\varphi}(s^{-1}x)=R_{\varphi}(x),~\forall x\in \mathcal{M}\overline{\otimes}\mathcal{N}.\]
Towards this end, it is enough to check this for elements of the form $\sum_{i=1}^na_i\otimes b_i$, where $a_i\in\mathcal{M}$ and $b_i\in\mathcal{N}$ for all $i=1,2,\ldots,n$. We now observe that
\begin{align*}
&\sum_{s\in\Gamma}\mu(s)sR_{\varphi}\left(s^{-1}\left(\sum_{i=1}^na_i\otimes b_i\right)\right)\\&=\sum_{s\in\Gamma}\mu(s)sR_{\varphi}\left(\left(\sum_{i=1}^ns^{-1}a_i\otimes s^{-1}b_i\right)\right)\\&=\sum_{s\in\Gamma}\sum_{i=1}^n\mu(s)sR_{\varphi}\left(s^{-1}a_i\otimes s^{-1}b_i\right)\\&=  \sum_{s\in\Gamma}\sum_{i=1}^n\mu(s)s\left(\varphi(s^{-1}a_i)s^{-1}b_i\right)\\&= \sum_{s\in\Gamma}\sum_{i=1}^n\mu(s)\varphi(s^{-1}a_i)s\left(s^{-1}b_i\right)\\&=\sum_{i=1}^n\left(\left(\sum_{s\in\Gamma}\mu(s)\varphi(s^{-1}a_i)\right)b_i\right)\\&=\sum_{i=1}^n\varphi(a_i)b_i=R_{\varphi}\left(\sum_{i=1}^na_i\otimes b_i\right).
\end{align*}
Given a $\mu$-stationary state $\psi$ on $\mathcal{N}$, a similar argument would work to show that the right slice map $R_{\psi}:\mathcal{M}\overline{\otimes}\mathcal{N}\to\mathcal{M}$  defined by $a\otimes b\mapsto a\psi(b)$ is $\mu$-stationary. We also note that a similar computation works in the minimal tensor product of $C^*$-algebras.
\end{example}
Let $\mathcal{A}$ be a separable $\Gamma$-$C^*$-algebra. Given a measurable $\Gamma$-space $(Y,\eta)$, let $\mathcal{C}_{Y,\mathcal{A}}$ denote the collection of $\eta$-almost everywhere measurable maps from $Y$ to the state space of $\mathcal{A}$, i.e.,
\[\mathcal{C}_{Y,\mathcal{A}}=\{\varphi:Y\to S(\mathcal{A}):\varphi\text{ is $\eta$-a.e measurable}\}.\]
There is a natural action of $\Gamma$ on $\mathcal{C}_{Y,\mathcal{A}}$ by left translations on both sides. In particular, for $s\in\Gamma$ and $\varphi\in \mathcal{C}_{Y,\mathcal{A}}$, $$(s.\varphi)(y)(a):=\varphi(s^{-1}y)(s^{-1}a),\  y\in Y_0, \ a\in\mathcal{A}.$$
Given a probability measure $\mu\in\text{Prob}(\Gamma)$, a fixed point in $\mathcal{C}_Y$ under the convolution of $\mu$ is called  $\mu$-stationary. That such a stationary map always exists is a standard Markov-Kakutani fixed point argument. 
\begin{proposition}
\label{proposition:almostexistence}
Let $\mathcal{A}$ be a separable $\Gamma$-$C^*$-algebra, $(Y,\eta)$ a measurable $\Gamma$-space, and $\mathcal{C}_{Y,\mathcal{A}}$ the associated collection of measurable maps. Given $\mu\in\text{Prob}(\Gamma)$, there exists $\varphi\in\mathcal{C}_{Y,\mathcal{A}}$ which is $\mu$-stationary in the sense that
\[\left(\mu*\varphi\right)(y)(a):=\sum_{s\in\Gamma}\mu(s)\varphi(s^{-1}y)(s^{-1}a)=\varphi(y)(a), \ y\in Y_0, \ a\in\mathcal{A}.\]
\begin{proof}
The set $\mathcal{C}_{Y,\mathcal{A}}$ is a convex set with respect to the pointwise addition induced from $S(\mathcal{A})$. In particular for $\varphi,\psi\in \mathcal{C}_{Y,\mathcal{A}}$ and $t\in (0,1)$, 
\[\left(t\varphi+(1-t)\psi\right)(y)(a)=t\varphi(y)(a)+(1-t)\psi(y)(a),\ y\in Y_0, \ a\in\mathcal{A}.\]
Since $\mathcal{A}$ is separable, there is an affine bijection between $\mathcal{C}_{Y,\mathcal{A}}$ and the convex set of unital completely positive maps $E:\mathcal{A}\to L^{\infty}(Y,\eta)$ (see \cite[Theorem~2.10]{bassi2020separable}). The collection of all unital completely positive maps $E:\mathcal{A}\to L^{\infty}(Y,\eta)$ is compact when endowed with the point ultraweak topology (see~\cite[Theorem~1.3.7]{BroOza08}). This induces on $\mathcal{C}_{Y,\mathcal{A}}$ a topology which makes it compact. To sum up, $\mathcal{C}_{Y,\mathcal{A}}$ is a compact convex $\Gamma$-space. Given $\mu\in\text{Prob}(\Gamma)$, consider the induced affine continuous map $M_{\mu}: \mathcal{C}_{Y,\mathcal{A}}\to \mathcal{C}_{Y,\mathcal{A}}$ defined by $\varphi\mapsto\mu*\varphi$. The claim now follows from the Markov-Kakutani fixed point theorem.
\end{proof}
\end{proposition}
Now, given a $\mu$-stationary element in $\mathcal{C}_{Y,\mathcal{A}}$, using the affine bijection between $\mathcal{C}_{Y,\mathcal{A}}$ and unital completely positive maps $E:\mathcal{A}\to L^{\infty}(Y,\eta)$, we obtain a $\mu$-stationary u.c.p map $E:\mathcal{A}\to L^{\infty}(Y,\eta)$.
In subsequent work, we will investigate further the existence and uniqueness of such $\mu$-stationary maps.
\begin{remark}
It follows from Proposition~\ref{proposition:almostexistence} that we can find a $\mu$-stationary $\eta$-a.e defined $\mu$-stationary measurable map $\varphi:Y\to \text{Prob}(X)$. Moreover, if we let $\nu=\text{bar}(\phi\eta)$, it is easy to see that $\nu=\int_Y\varphi(y)d\eta(y)$. We do not know if $\varphi$ is the disintegration map associated with some factor map $\pi: (X,\nu)\to (Y,\eta)$.
\end{remark}
\section{Relative proximality}
\label{sec:relprox}
The Poisson boundary is probably the most famous example of a stationary dynamical system. The Poisson boundary and its $\G$-factors, so-called $\mu$-boundaries, come with another important property, namely \textit{proximaility}. In the present section, we generalise this concept to $\G$-factor maps. Again, it coincides with the usual definition of a proximal system if the factor map goes to a one-point space.

To this end, let us denote \begin{equation}\label{eq: lim}\nu_{y,\omega}:=\lim\limits_{n\to\infty} g_1\ldots g_n\nu_{g_n^{-1}\ldots g_1^{-1} y}    
\end{equation} 
for $\omega=(g_1,g_2,\ldots)\in\G^{\mathbb{N}}$ and $y \in Y$, provided this limit exists in the weak*-topology on $\prob(X)$, the space of Borel probability measures on $X$.

The above disintegration limits are in some sense equivariant with respect to a ``two-sided'' action as follows:
\begin{lemma}[Equivariance]\label{lem: equiv} Let $X$ be a compact Hausdorff space with a continuous $G$-action.
Let $\pi:(X,\nu)\longrightarrow (Y,\eta)$  be a continuous $\G$-factor map.
    Fix $y\in Y$ such that $\nu_y$ exists. Then
$$
        \g\nu_{\g^{-1} y, (\g_1,\g_2,\ldots)}=\nu_{y, (\g,\g_1,g_2,\ldots)}
    $$ for all $\g\in \G$ and  $\mu^{\mathbb{N}}$-a.e. $(\g_1,\g_2,\ldots) \in\G^{\mathbb{N}}$. 
\end{lemma}
\begin{proof}
By continuity of the action, we obtain
$$ \g\nu_{\g^{-1} y, (\g_1,\g_2,\ldots)} =
\lim\limits_{n\to\infty} \g\g_1\ldots\g_n\nu_{\g_n^{-1}\ldots\g_1^{-1}\g^{-1}y}=\nu_{y, (\g,\g_1,g_2,\ldots)},
$$ which proves the claim.
\end{proof}
A prior, it is not obvious that the limit in equation~(\ref{eq: lim}) exists. We show below that they do and, thereby, build the connection between relative stationary systems and the theory of random walks. 
\begin{theorem}[Existence and barycenter]\label{thm: lim ex}
Let $X$ be a compact Hausdorff space with a continuous $G$-action. Let $\mu$ be a non-degenerate and admissible Borel probability measure on $G$.
Let $\pi:(X,\nu)\longrightarrow (Y,\eta)$ be relatively $\mu$-stationary and continuous. Let $\Omega:=\G^{\mathbb{N}}$.
    Fix $y\in Y$ such that $\nu_y$ exists. Then there exists a co-null set $\Omega_y\in\mathcal{B}(\G)^{\mathbb{N}}$ such that the limit $\nu_{y,\omega}$ exists for all $\omega\in \Omega_y$. Moreover, the barycenter equation holds:
    \[\int_{\Omega}\nu_{y,\omega}\, d\mu^{\mathbb{N}}(\omega)=\nu_y.\]
\end{theorem}
\begin{proof}
We shall prove it in a series of steps.\\ 
We denote by $\text{Har}(G,\mu)$, the space of bounded right $\mu$-harmonic functions i.e., $\varphi: G \to \mathbb{R}$ is $\mu$-harmonic if 
\[\varphi(g)=\int_G \varphi(gg')\, d\mu(g'), g \in G.\]
There is a natural action of $G$ on $\text{Har}(G,\mu)$ namely, $g_1\varphi(g)=\varphi(g_1^{-1}g)$ for $g,g_1 \in G$. 
\\ 
\textit{Step-1: 
Let $\varphi \in \text{Har}(G,\mu)$. Then, for $\mu^{\mathbb{N}}$-a.e $\omega \in \Omega$, the limit $$\lim_n \varphi(\omega^n)=\Tilde{\varphi}(\omega)$$ exists and belongs to $L^{\infty}(\Omega,\mu^{\mathbb{N}})$. Moreover,
\[\varphi(e)=\int_{\Omega}\Tilde{\varphi}(\omega)\, d\mu^{\mathbb{N}}(\omega).\]}
\noindent\\
The above claim is well-known, but we include its proof for completeness.
For each $n\in\mathbb{N}$, let $X_n:\Omega\to G$ denote the projection map onto the $n$\textsuperscript{th}-coordinate, and let $X_0=e$. Let $\mathcal{F}_0=\sigma(\{e\})$, $\mathcal{F}_n=\sigma\left(X_0,X_1,\ldots,X_n\right)$. Then, it is clear that $\mathcal{F}_n \subseteq \mathcal{F}_{n+1}$ for all $n \ge 0$. For each $\omega\in\Omega$, let $\omega^n:=\omega_1\omega_2\ldots\omega_n$, where $\omega=(\omega_1,\omega_2,\ldots)$. Now, for $\varphi \in \text{Har}(G,\mu)$, we define $\varphi_n: \Omega \to \mathbb{R}$ by $\varphi_n(\omega)=\varphi(w^n)=\varphi(\omega_1\omega_2\ldots\omega_n)$ (the projection into the first n-coordinates), $\varphi_0=e$. Since $\varphi$ is bounded, $\varphi_n$ is also bounded. 
All we need to show now is that $\mathbb{E}\left(\varphi_{n+1}|\mathcal{F}_n\right)=\varphi_n$. Once we establish this, by the Martingale convergence Theorem, the limit $\lim_n\varphi_n(\omega)$ will exist $\mu^{\mathbb{N}}$-a.e. and we denote this by $\tilde{\varphi}(\omega)$. This limit will be in $L^{\infty}(\Omega, \mu^{\mathbb{N}})$ since $\varphi$ is itself bounded. Moreover, 
\[\int_{\Omega}\varphi_n(\omega)\, d\mu^{\mathbb{N}}(\omega)=\int_{\Omega}\varphi_o(\omega)\, d\mu^{\mathbb{N}}(\omega)=\varphi(e).\]
Since $\varphi_n(\omega)\to \Tilde{\varphi}(\omega)$ pointwise $\mu^{\mathbb{N}}$-a.e, using the dominated convergence theorem, we will obtain that \[\int_{\Omega}\Tilde{\varphi}(\omega)\, d\mu^{\mathbb{N}}(\omega)=\varphi(e).\]
We now claim that $\mathbb{E}\left(\varphi_{n+1}|\mathcal{F}_n\right)=\varphi_n$. To establish this, we need to show that for all $F \in \mathcal{F}_n$, $$\int_F \varphi_{n+1}(\omega)d\mu^{\mathbb{N}}(\omega)=\int_F \varphi_n(\omega)d\mu^{\mathbb{N}}(\omega)$$ holds. Then, by uniqueness of $\mathbb{E}\left(\varphi_{n+1}|\mathcal{F}_n\right)$, we will obtain the claim. Towards that end, we can look at the basic generating sets of the form $F= {A_1}\times A_2\times \ldots \times A_n\times G\times\ldots$ which is in $\mathcal{F}_n$. 
By definition, 
\begin{align*}
    \int_F\varphi_{n+1}(\omega)d\mu^{\mathbb{N}}(\omega)&=\int_F \varphi(\omega_1\omega_2\ldots\omega_{n+1})\ d\mu^{\mathbb{N}}(\omega).
\end{align*}
Now,
\begin{align*}
    \int_F \varphi(\omega_1\omega_2\ldots\omega_{n+1})d\mu^{\mathbb{N}}(\omega)&=\int_{A_1} \ldots \int_{A_n}d\mu\int_{G}\varphi(\omega_1\ldots\omega_{n+1})\ d\mu(\omega_{n+1})\ d\mu(\omega_n)\ d\mu(\omega_1).
\end{align*}

Using harmonicity, we see that
\begin{align*}
    \int_{G}\varphi(\omega_1\omega_1\ldots\omega_{n+1})\, d\mu(\omega_{n+1})&=\varphi(\omega_1\omega_2\ldots\omega_n).
\end{align*}
Hence, we get that 
\begin{align*}
\int_F \varphi(\omega_1\omega_2\ldots\omega_{n+1})d\mu^{\mathbb{N}}(\omega)&=\int_{A_1} \ldots \int_{A_n} \int_{G}\varphi(\omega_1\ldots\omega_{n+1})\, d\mu(\omega_{n+1}) \, d\mu(\omega_{n})\, \ldots \, d\mu(\omega_1)\\&=\int_{A_1\times A_2\times\ldots\times A_n}\varphi(\omega_1\omega_2\ldots\omega_n)\, d\mu^{\mathbb{N}}(\omega)
\\&=\int_F \varphi_n(\omega)\, d\mu^{\mathbb{N}}(\omega).
\end{align*}
This proves the first step.\\
\textit{Step-2: Let $y\in Y$ be such that $\nu_y$ exists. Consider the map $P_{\nu_y}: L^{\infty}(\pi^{-1}(\{y\}),\nu_y) \to L^{\infty}(G,\mu)$ defined by
\[P_{\nu_y}(f)(g):=\int_X f(gx)\, d\nu_{\g^{-1}y}(x),\ g\in G.\]
Then, for every $f\in L^{\infty}(\pi^{-1}(\{y\}),\nu_y)$, $P_{\nu_y}(f)\in \text{Har}(G,\mu)$.}\\
The map $P_{\nu_y}$ is a natural generalisation of the Furstenberg transform and was already considered by \cite[Proposition 4.3.1]{creutz2016contractive}.
First note that $\nu_y$ and $g\nu_{g^{-1}y}$ share the same nullsets for all $g\in G$, because $\nu_y=\int_G g\nu_{g^{-1}y} \, d\mu(g)$ and $\mu$ is non-degenerate. Therefore, the map $P_{\nu_y}$ is well defined.
For every $f \in L^{\infty}(\pi^{-1}(\{y\}),\nu_y)$, we observe that
\begin{align*}
P_{\nu_y}(f)(g)&=\int_X f(gx)d\nu_{\g^{-1}y}(x)
%\\&=\int_X f(x)dg.\nu_y(x)
\\&=\int_X f(x)dg\nu_{g^{-1}y}(x)\\&=g\nu_{g^{-1}y}(f).    
\end{align*}
Now, 
\begin{align*}
\int_G  P_{\nu_y}(f)(gh)\, d\mu(h)&= \int_G\int_X f(x)\, dgh\nu_{h^{-1}g^{-1}y}(x)\, d\mu(h)\\&= \int_G h\nu_{h^{-1}(g^{-1}y)}(g^{-1}f)\, d\mu(h)\\&=\nu_{g^{-1}y}(g^{-1}f)\\&=g\nu_{g^{-1}y}(f)=P_{\nu_y}(f)(g).   
\end{align*}
Therefore, $P_{\nu_y}(f)\in \text{Har}(G,\mu)$ for $\eta$-a.e $y\in Y$.\\
\textit{Step-3: For $\mu^{\mathbb{N}}$-a.e-$\omega \in \Omega$, the limit $\lim_n\omega^n\nu_y=\nu_{y,\omega}$ exists in the weak$^*$-topology of $P(X)$. Moreover, 
\[\nu_y=\int_{\Omega}\nu_{y,\omega}d\mu^{\mathbb{N}}(\omega).\]}
It follows from $\textit{Step-2}$ that $P_{\nu_y}(f)\in \text{Har}(G,\mu)$ for $\eta$-a.e $y\in Y$ and for every $f\in L^{\infty}(\pi^{-1}(\{y\}),\nu_y)$. By passing to a compact model, we can assume that $X$ is metrizable; hence, $C(X)$ is separable. Choose a countable dense subset $C_0 \subset C(X)$. For each $f \in C_0$, since $P_{\nu_y}(f)\in \text{Har}(G,\mu)$, it follows from \textit{Step-1} that we can find a $\mu^{\mathbb{N}}$-conull set $\Omega_f$ such that \[\lim_nP_{\nu_y}(f)(\omega^n)\text{ exists for all } \omega \in \Omega_f.\]
Let $\Omega_y=\cap_{f \in C_0}\Omega_f$. Observe that $\mu^{\mathbb{N}}(\Omega_y)=1$. For $\omega\in \Omega$, and $\omega^n=\omega_1\omega_2\ldots\omega_n$, we denote $\left(\omega_1\omega_2\ldots\omega_n\right)^{-1}$ by $\omega^{-n}$. Now, for all $\omega \in \Omega_y$ and for $f \in C_0$, the following limit exists:
\[\lim_nP_{\nu_y}(f)(\omega^n)=\lim_n\int_Xf(\omega^nx)d\nu_{\omega^{-n}y}(x)=\lim_n\int_Xf(x)d(\omega^n\nu_{\omega^{-n}y})(x).\]
For each $\omega \in \Omega_0$, we define a linear functional $\psi_{\omega,y}: C_0 \to \mathbb{C}$ defined by
\[\psi_{\omega,y}(f)=\lim_n\int_Xf(x)d(\omega^n\nu_{\omega^{-n}y})(x).\]
Since $\psi_{\omega,y}$ is bounded, it extends to a bounded linear functional $\tilde{\psi}_{\omega,y}$ on $C(X)$. Using Riesz-Representation theorem, we obtain a probability measure $\nu_{y,\omega}$ on $X$ such that
\[\tilde{\psi}_{\omega,y}(f)=\int_Y fd\nu_{y,\omega},\  f\in C(X), \ \omega \in \Omega_0.\]
Moreover, 
\begin{align*}
\int_{\Omega}d\nu_{y,\omega}(f)d\mu^{\mathbb{N}}(\omega)&=\int_{\Omega}\lim_n\langle f,\omega^n\nu_{\omega^{-n}y}\rangle d\mu^{\mathbb{N}}(\omega)\\&=\lim_n\int_{\Omega}\langle f,\omega^n\nu_{\omega^{-n}y}\rangle d\mu^{\mathbb{N}}(\omega)\\&=\lim_n\int_G\int_{\{\Omega: \omega^n=g\}}\langle f,\omega^n\nu_{\omega^{-n}y}\rangle d\mu^{\mathbb{N}}(\omega)d\mu(g)\\&=\lim_n\int_G\langle f,g\nu_{g^{-1}y}\rangle\mu^{\mathbb{N}}\left(\{\Omega: \omega^n=g\}\right)d\mu(g)\\&=\lim_n\left\langle f,\int_Gg\nu_{g^{-1}y}d\mu^{*n}(g)\right\rangle\\&=\lim_n\left\langle f,\nu_y\right\rangle=\nu_y(f).
    \end{align*}

\end{proof}
We investigate the conull sets $\Omega_y$ from the previous proof and put enough of these sets together to establish that the limit in equation~\eqref{eq: lim} exists almost everywhere.
\begin{lemma}\label{lem: conull Omega}
Let $Y_0$ be a conull set such that $\nu_y$ exists for all $y\in Y_0$. For $y\in Y_0$ let $\Omega_y$ be a conull sets in $\G^{\mathbb{N}}$ such that $\nu_{y,\omega}$ exists for all $\omega\in \Omega_y$. Then $\bigcup_{y\in Y_0}\{y\}\times \Omega_y$ is measurable and $\eta\otimes \mu^{\mathbb{N}}(\bigcup_{y\in Y_0}\{y\}\times \Omega_y)=1 .$
\end{lemma}

\begin{proof} W.l.o.g. we may assume that $X$ is compact metrizable by Proposition~\ref{compactmodelsforfactormaps}.
By Proposition~\ref{proposition:measurabilityofthedisint}, the almost everywhere defined map $y\mapsto \nu_y$ is measurable. Let $n\in\mathbb{N}$. Define $f_n:Y\times\Omega\to\text{Prob}(X)$ by $f_n(y,\omega)=\omega^{n}\nu_{\omega^{-n}y}$. Since $X$ is metrizable, $\text{Prob}(X)$ is a separable metric space. Since $f_n$'s are measurable, using \cite[Lemma~1.11]{Kallenberg}, the set $\{(y,\omega): f_n(y,\omega)\text{ converges}\}$ is measurable. Now, $\bigcup_{y\in Y_0}\{y\}\times \Omega_y$ is precisely the domain where the functions $f_n$ converge and thus this set is measurable. Moreover, 
$$\eta\otimes \mu^{\mathbb{N}}\Big(\bigcup_{y\in Y_0}\{y\}\times \Omega_y\Big)=\int_{Y} \delta_{\tilde{y}}\otimes  \mu^{\mathbb{N}}  \Big(\bigcup_{y\in Y_0}\{y\}\times \Omega_y\Big) \ d\eta(\tilde{y})
$$ $$=\int_{Y}   \mu^{\mathbb{N}}  (\Omega_{\tilde{y}}) \ d\eta(\tilde{y})=1.
$$
Here we used that  $\delta_{\tilde{y}}\otimes  \mu^{\mathbb{N}}  (\bigcup_{y\neq \tilde{y}}\{y\}\times \Omega_y) \leq \delta_{\tilde{y}}\otimes  \mu^{\mathbb{N}}  (\bigcup_{y\neq \tilde{y}}\{y\}\times \G^{\mathbb{N}})=0.$    
\end{proof}

With the collection of our previous results, we can now define a relative boundary map as follows.
\begin{corollary}\label{corollary: bnd to prob}
Let $\pi: (X,\nu)\longrightarrow (Y,\eta)$ be relatively $\mu$-stationary map. Then there exists an almost everywhere defined map $$ \widetilde{\bnd}: Y\times \G^{\mathbb{N}} \longrightarrow \prob(X)$$  which is $\G$-equivariant and  measurable.
\end{corollary}

\begin{proof} By Proposition~\ref{compactmodelsforfactormaps} we may pass to compact models, i.e. $\pi$ is continuous, $X$ is compact Hausdorff and $G$ acts continuously on $X$.
    By Theorem~\ref{thm: lim ex},   we obtain a map $$\bigcup_{y\in Y_0}\{y\}\times \Omega_y \longrightarrow \prob(X),
    $$ 
    $$ (y,\omega)\mapsto \nu_{y,\omega},$$
    where $\Omega_y$ is a co-null set of $\Omega:=\G^{\mathbb{N}}$ given by Theorem~\ref{thm: lim ex} and $Y_0$ is a conull set such that $\nu_y$ exists for every $y\in Y_0$. By Lemma~\ref{lem: equiv} this map is $\G$-equivariant, since 
    $$\nu_{\g y,\g \omega}=  \g \nu_{y,\omega}.
    $$ Due to Proposition~\ref{proposition:measurabilityofthedisint} and Lemma~\ref{lem: conull Omega} the map is almost everywhere defined and measurable.
\end{proof}

Finally, we are in a position to define a relative version of proximality.

\begin{definition}[Relative proximal]\label{def: prox}
    We call  $\pi:(X,\nu)\longrightarrow (Y,\eta)$ relatively proximal if for $\eta$-a.e. $y\in Y$ and  for $\mu^{\mathbb{N}}$-a.e. $ \omega\in \G^{\mathbb{N}}$, the limit $\nu_{y,\omega}$ exists and 
    $$\nu_{y,\omega}=\delta_{\bnd(y,\omega)},$$ where $\bnd(y,\omega)$ is an element in $X$ which depends on $y$ and $\omega$.
\end{definition}

A related notion has been studied by \cite{creutz2016contractive}
under the name of relatively contractive maps. Contrary to our definition, he assumes that there is one sequence that shrinks the measure to a Dirac measure.

Proximaility implies that the map $\widetilde{\bnd}: Y\times \G^\mathbb{N}\longrightarrow \prob(X)$ in fact goes to $X$ rather than $\prob(X)$.
\begin{corollary}\label{corollary: bnd map}
Let $\pi: (X,\nu)\longrightarrow (Y,\eta)$ be relatively proximal. Then there exists a $\G$-equivariant, measurable almost everywhere defined map $$ \bnd: Y\times \G^\mathbb{N}\longrightarrow X$$  which  satisfies  $\pi\circ \bnd= \pr_1$, where $\pr_1: Y\times \G^\mathbb{N}\longrightarrow Y$ denotes the projection.
\end{corollary}
\begin{proof}
By Corollary~\ref{corollary: bnd to prob} we obtain a $\G$-equivariant, measurable map to $\prob(X)$, given by $(y,\omega)\mapsto \nu_{y,\omega}$ for $\mu^{\mathbb{N}}$-a.e. $\omega\in \G^\mathbb{N} $ and $\eta$-a.e. $y\in Y$. Now by assumption, we have $\nu_{y,\omega}=\delta_{\bnd(y,\omega)}$ almost everywhere. Thus, identifying $\delta_{\bnd(y,\omega)}$ with $\bnd(y,\omega)$ we obtain the desired map to $X$.
Since $\supp(\nu_{y,\omega})\subseteq \pi^{-1}(y)$, we obtain $\bnd(y,\omega)\subseteq \pi^{-1}(y)$ and therefore $ \pi(\bnd(y,\omega))=y$ for almost all $y\in Y$, almost all $\omega\in\G^{\mathbb{N}}$.
\end{proof}

\subsection{Heredity to intermediate factors}
Similar to relative $\mu$-stationarity, the property of relative proximality passes to intermediate factors (for a similar result in Creutz's notion, see \cite[Theorem 4.7]{creutz2016contractive}).
\begin{proposition}[Intermediate factors]\label{proposition: IF prox}
       Let $\pi: (X,\nu)\longrightarrow (Y,\eta)$ be relatively proximal and let $$         \xymatrix{ (X,\nu)\ar[dd]^{\pi}\ar[dr]^{\kappa} & 
        \\ 
        & (Z,\zeta)\ar[dl]^{\sigma}\\
        (Y,\eta)
        }        $$ 
               be a commutative diagram of $\G$-factor maps. Then, $\sigma$ is relatively proximal.
\end{proposition}

\begin{proof} Without any loss of generality, we may assume that all the maps are continuous. 
First note that $(\kappa \nu)_y=\kappa(\nu_y)$ for almost every $y\in Y$, since $\int_Y \kappa(\nu_y)\, d\eta(y)= \kappa \nu
$.
Thus, by continuity $$\zeta_{y,\omega}=(\kappa \nu)_{y,\omega}=\kappa (\nu_{y,\omega})=\kappa \delta_{\bnd(\omega,y)}$$
for almost all $y\in Y$ and $\omega\in \Omega$. This means that $\sigma$ is relatively proximal.    
\end{proof}

\subsection{Only trivial endomorphisms}
It is known that proximal systems cannot have non-trivial factor maps onto themselves. We show that this also holds true of our relative notion.

\begin{proposition}\label{proposition: id map}
Let $\pi: (X,\nu)\longrightarrow (Y,\eta)$ be a relatively proximal $\G$-factor map. Assume there exists a $\G$-factor map $\alpha: (X,\nu)\longrightarrow (X,\nu)$ such that $\pi\circ \alpha=\pi$. Then $\alpha$ is the identity map almost everywhere.
   \end{proposition}

\begin{proof}
We modify the arguments in \cite[Proposition 3.2]{GF09} and adapt them to our purpose. Let $$\bary: \prob(\prob(X))\longrightarrow \prob(X)$$ denote the barycenter map. We refer the reader to \cite[Chapter-II]{glasner1976proximal} for the definition of the barycenter map. For every $f\in C(X)$, let $\Phi_f:\text{Prob}(X)\to\mathbb{C}$ be defined by $$\Phi_f(\xi) = \int f \, d\xi,~\xi\in \prob(X).$$  By the properties of the barycenter map, we know that 
\begin{equation}\label{eq: Phi}
    \Phi_f(\bary( \Theta))= \int \Phi_f\, d\Theta
\end{equation} for $\Theta\in \prob(\prob(X))$
(see for instance \cite[Proposition 2.1]{glasner1976proximal}).
Let $$\iota: x\mapsto \frac{1}{2}\Big(\delta_{x}+\delta_{\alpha (x)}\Big)$$ be a map from $X$ to $\prob(X)$. Then we have $$
    \bary(\iota \nu)=\nu.$$
Indeed, by (\ref{eq: Phi}) we see that for all $f\in C(X)$,
$$
\bary(\iota\nu)(f)=\Phi_f(\bary(\iota\nu))=
\int \Phi_f \, d \iota\nu
$$
$$= \int \int f\, d \xi \, d \iota\nu(\xi)=
\int \frac{1}{2}\Big( f(x) + f(\alpha(x))\Big) \, d \nu(x)
$$
$$
= \frac{1}{2}\Big(\int f(x) \, d \nu(x) +\int  f(x) \, d \alpha\nu(x) \Big)
$$
$$= \nu(f),
$$
since $\alpha\nu=\nu$.
Moreover, we claim that
\begin{equation}\label{eq: iota dis}\bary(\iota(\nu_y))=\nu_{y}\end{equation} for $\eta$-almost every $y\in Y$. We first show that $\supp(\bary(\iota(\nu_y))
\subset\pi^{-1}(y)$. Let $f\in C(X)$ be such that $\supp(f)\cap \pi^{-1}(y)=\emptyset$. Since $\pi\circ\alpha=\pi$, we see that $\pi^{-1}(y)=\alpha^{-1}(\pi^{-1}(y))$ for $y\in Y$. Therefore, we obtain that $$\emptyset=\alpha^{-1}(\pi^{-1}(y))\cap \alpha^{-1}(\supp(f))=\pi^{-1}(y)\cap\alpha^{-1}(\supp(f)).$$ Now, it follows that \begin{align*}
 \bary(\iota\nu_y)(f)&=
\int_{\text{Prob}(X)} \Phi_f(\xi) \, d \iota\nu_y(\xi)\\&= \int_X\Phi_f\circ \iota(x)d\nu_y(x)\\&=
\int_X \frac{1}{2}\Big( f(x) + f(\alpha(x))\Big) \, d \nu_y(x)
\\&= \frac{1}{2}\int_X f(x) \, d \nu_y(x) +\frac{1}{2}\int_X  f(\alpha(x)) \, d \nu_y(x)\\&=\frac{1}{2}\int_{\pi^{-1}(y)} f(x) \, d \nu_y(x) +\frac{1}{2}\int_{\pi^{-1}(y)} f(\alpha(x)) \, d \nu_y(x)\\&=\frac{1}{2}\int_{\pi^{-1}(y)} f(\alpha(x)) \, d \nu_y(x)&(\text{$\supp(f)\cap \pi^{-1}(y)=\emptyset$})\\&=0 &(\text{$\alpha^{-1}(\supp(f))\cap \pi^{-1}(y)=\emptyset$})  
.\end{align*} This shows that $\supp(\bary(\iota(\nu_y))
\subset\pi^{-1}(y)$.
The claim now follows by the uniqueness of the disintegration measures and 
$$
\int_Y \bary(\iota(\nu_y))(f)\, d\eta(y)= 
\int_Y \Phi_f(\bary(\iota(\nu_y)))\, d\eta(y)
$$
$$=
\int_Y \int \Phi_f\circ \iota\, d \nu_y \, d\eta(y)
= \int \Phi_f\circ \iota\, d \nu =\nu(f).
$$
Now, by assumption $\pi$ is relatively proximal, hence for almost every $y\in Y$ and almost every $\omega\in \G^{\mathbb{N}} $, we have 
$$\lim\limits_{n\to \infty} \omega^n\nu_{\omega^{-n}y} = \delta_{\bnd(y,\omega)}
$$ for $\bnd(y,\omega)\in X$,
where $\omega^n:=\g_1\cdots \g_n$ and $\omega^{-n}:=(\omega^{n})^{-1}$, for $\omega=(\g_1,\g_2,\ldots)$.
Hence, 
$$
\bary(\iota \lim\limits_{n\to \infty} \omega^n\nu_{\omega^{-n}y})=\iota(\bnd(y,\omega))=\frac{1}{2}\Big(\delta_{\bnd(y,\omega)}+\delta_{\alpha(\bnd(y,\omega))}\Big).
$$
On the other hand, since the barycenter map and $\iota$ are  $\G$-equivariant and weak*-continuous, we obtain  that 
\begin{align*}
\bary(\iota \lim\limits_{n\to \infty} \omega^n\nu_{\omega^{-n}y})&=
 \lim\limits_{n\to \infty} \omega^n \bary(\iota(\nu_{\omega^{-n}y}))\\&\stackrel{\text{Eq. }\eqref{eq: iota dis}}{=}
  \lim\limits_{n\to \infty} \omega^n \nu_{\omega^{-n}y} \\&=\delta_{\bnd(y,\omega)}.
\end{align*}
Combining these two series of equations, we see that $$
\delta_{\bnd(y,\omega)}=\frac{1}{2}\Big(\delta_{\bnd(y,\omega)}+\delta_{\alpha(\bnd(y,\omega))}\Big).
$$ Thus $\bnd(y,\omega)=\alpha(\bnd(y,\omega))$ for almost every $y\in Y$ and almost every $\omega\in G^{\mathbb{N}}$.
We observe that $\bnd$ is surjective up to a null set since  $$\nu=\int_{G^\mathbb{N}}\int_Y \nu_{y,\omega}\, d\eta(y) \, d\mu^{\mathbb{N}}(\omega)=\int_{G^\mathbb{N}}\int_Y \delta_{\bnd(y,\omega)}\, d\eta(y) \, d\mu^{\mathbb{N}}(\omega),$$ Hence, $X=\bnd(Y\times G^\mathbb{N})$ up to a null set, thereby implying that $\alpha=\id$ on $X$ up to a null set.
\end{proof}

\section{Relative \texorpdfstring{$\mu$}{mu}-boundaries}
\label{sec:relboundary}
Compared to topological ones, a key benefit of measurable boundaries is their ease of concrete identification. Indeed, the past few decades have seen significant progress, culminating in the concrete realisation of the Poisson boundary for numerous groups that naturally emerge as symmetries of geometric shapes.
A $\mu$-boundary is a $\mu$-stationary and proximal dynamical system. We extend this definition to maps between dynamical systems.

\begin{definition}[Relative $\mu$-boundaries]\label{def: rel bnd}
    A relatively proximal and relatively $\mu$-stationary $\G$-factor map $\pi:(X,\nu)\longrightarrow (Y,\eta)$ is called a \textit{relative $\mu$-boundary}.
\end{definition}
We give an example (see \cite[Theorem~3.4]{Naghavi}) of a product space where the action is not a product action and yet is a relative boundary. To this end, we restrict our considerations to a countable group $\Gamma$.
Given a $\Gamma$-space $X$ and a subgroup $\Lambda\le \Gamma$, recall that a cocycle of the action in $\Lambda$ is a map $\alpha:\Gamma\times X\to \Lambda$ such that 
\[\alpha(\gamma_1\gamma_2,x)=\alpha(\gamma_1,\gamma_2x)\alpha(\gamma_2,x), \ \forall \gamma_1,\gamma_2\in\Gamma,\ x\in X.\]
\begin{example} 
 Let $\Lambda\le \Gamma$ be a finite index subgroup. Let $Y$ be a $\Lambda$-space. Let $T=\{\Lambda,t_1\Lambda,\ldots,t_n\Lambda\}$ be a representative for $\Gamma/\Lambda$. Define $\alpha:\Gamma\times\Gamma/\Lambda\to \Lambda$ by 
 \[\alpha(\gamma, t_i\Lambda)=\lambda,\]
 where $\lambda$ is uniquely determined by $\gamma t_i\lambda\in T$. 
 Let $X=\Gamma/\Lambda\times Y$. Define the action $\Gamma\curvearrowright X$ by
\[\gamma(t_i\Lambda, y)=\left(\gamma t_i\alpha(\gamma,t_i\Lambda)\Lambda,\alpha(\gamma,t_i\Lambda)^{-1}y \right), \ y\in Y, \ \gamma\in\Gamma.\]
The space $X$ with the above-defined action is the induced $\Gamma$-space of the $\Lambda$-space $Y$. Consider now the projection map $\pr_{\Gamma/\Lambda}:\Gamma/\Lambda\times Y\to \Gamma/\Lambda$ defined by $(t_i\Lambda, y)\mapsto t_i\Lambda$ for each $y\in Y$. We now show that $\pr_{\Gamma/\Lambda}$ is $\Gamma$-equivariant with respect to the action defined above. We first observe that
\[\pr_{\Gamma/\Lambda}\left(\gamma(t_i\Lambda, y)\right)=\gamma t_i\alpha(\gamma,t_i\Lambda)\Lambda.\]
Letting $\Lambda_i=t_i\Lambda t_i^{-1}$, we observe that $\Lambda_i$ is a finite index subgroup of $\Lambda$ with the representatives as $\{t_i^{-1}\Lambda_i, t_1t_i^{-1}\Lambda_i,\ldots,t_nt_i^{-1}\Lambda_i\}$ for $\Gamma/\Lambda_i$. In particular, $\Gamma=\bigsqcup_{j=1}^nt_jt_i^{-1}\Lambda$. We can write $\gamma=t_jt_i^{-1}\lambda_j$ for some $\lambda_j\in\Lambda$. Therefore, $\gamma t_i=t_jt_i^{-1}\lambda_jt_i\in t_j\Lambda$ (since $t_i^{-1}\lambda_jt_i\in\Lambda$). Hence, $\gamma t_i\Lambda=t_j\Lambda$. On the other hand, $\alpha(\gamma, t_i\Lambda)=e$. Consequently, it follows that
\[\gamma \pr_{\Gamma/\Lambda}\left(t_i\Lambda, y\right)=\gamma t_i\Lambda=\gamma t_i\alpha(\gamma,t_i\Lambda)\Lambda=\pr_{\Gamma/\Lambda}\left(\gamma(t_i\Lambda, y)\right).\]
Let $\mu\in\text{Prob}(\Lambda)$ and $\nu\in\text{Prob}(Y)$ be a $\mu$-stationary measure. Let $\eta\in\text{Prob}(\Gamma/\Lambda)$. Then $(\eta\times\nu)_{t_i\Lambda}=\delta_{t_i\Lambda}\times\nu$. Recall that $\Lambda_i=t_i\Lambda t_i^{-1}$. Let $\tilde{\mu}\in \text{Prob}(\Lambda_i)$ be defined by $\tilde{\mu}(t_is t_i^{-1})=\mu(s)$.
Since $ \alpha(t_ist_i^{-1}, t_i\Lambda)=s^{-1}$, we observe that 
\begin{equation}
\label{eq:actiondefined}
t_ist_i^{-1}.(t_i\Lambda, y)=\left(t_is\alpha(t_ist_i^{-1}, t_i\Lambda)\Lambda, \alpha(t_ist_i^{-1}, t_i\Lambda)^{-1}y\right)=\left(t_i\Lambda, sy\right).\end{equation}
We also note that
\[t_is^{-1}t_i^{-1}.(t_i\Lambda)=t_is^{-1}\alpha(t_is^{-1}t_i^{-1}, t_i\Lambda)\Lambda=t_i\Lambda.\]
Therefore,
\begin{align*}
&\int_{\Gamma/\Lambda\times Y}h(t\Lambda, y)d(t_ist_i^{-1}(\delta_{t_i\Lambda}\times\nu))\\&= \int_{\Gamma/\Lambda\times Y}h(t_ist_i^{-1}.(t\Lambda, y))d\left(\delta_{t_i\Lambda}\times \nu\right)\\&=\int_Y  h(t_ist_i^{-1}.(t_i\Lambda, y))d\nu\\&=\int_Y h(t_i\Lambda, sy)d\nu\\&=\int_Y h(t_i\Lambda, y)d(s\nu).
\end{align*}
Let us rewrite $h_{t_i\Lambda}\in L^{\infty}(Y,\eta)$ as $h_{t_i\Lambda}(y)=h(t_i\Lambda,y)$. Hence, it follows from above that
\begin{align*}
\sum_s \tilde{\mu}(t_ist_i^{-1})t_ist_i^{-1}\left(\delta_{t_i\Lambda}\times\nu\right)(h)&=\sum_s\mu(s)s\nu(h_{t_i\Lambda})\\&=\nu(h_{t_i\Lambda})\\&=\left(\delta_{t_i\Lambda}\times\nu\right)(h)  . 
\end{align*}
Consequently, it follows that $(\eta\times\nu)_{t_i\Lambda}$ is $\tilde{\mu}$-stationary. Let us further assume that $(Y,\nu)$ is a $\mu$-boundary and $\tilde{\mu}$ as above. Then, for $\omega\in \Lambda_i^{\mathbb{N}}$, we see that $\omega_n=t_is_1s_2\ldots s_n t_i^{-1}$. Therefore, arguing similarly to the above, we see that 
\[\int_{\Gamma/\Lambda\times Y}h(t\Lambda, y)d(\omega_n(\delta_{t_i\Lambda}\times\nu))=\int_Yh_{t_i\Lambda}(y)d(s_1s_2\ldots s_n\nu).\]
Since $\nu$ is a $\mu$-boundary, it follows that $\lim_n s_1s_2\ldots s_n\nu=\delta_{\text{bnd}(\Tilde{\omega})}$, where $$\Tilde{\omega}=(s_1,s_2, \ldots,s_n,\ldots)\in \Lambda^{\mathbb{N}},$$ and $\text{bnd}:\bigotimes(\Lambda, \mu)\to \text{Prob}(Y)$ is the associated boundary map. Therefore, using equation~\eqref{eq:actiondefined}, we see that
\[\lim_n\omega^n.((\delta_{t_i\Lambda}\times\nu))=\omega^n(\delta_{t_i\Lambda}\times\nu)=\delta_{t_i\Lambda}\times \delta_{\text{bnd}(\Tilde{\omega})}.\] 
This shows that $\pr_{\Gamma/\Lambda}:\Gamma/\Lambda\times Y\to \Gamma/\Lambda$ is a relative $\tilde{\mu}$-boundary. 
\end{example}

\subsection*{Heredity to intermediate factors}
Unsurprisingly, every intermediate factor of a relative $\mu$-boundary is a $\mu$-boundary. We obtain the following by combining Proposition~\ref{proposition: IF} along with \ref{proposition: IF prox}.
\begin{corollary}\label{cor: ift is rel prox}
Let $\pi:(X,\nu)\longrightarrow (Y,\eta)$ be a relative $\mu$-boundary and assume that $$         \xymatrix{ (X,\nu)\ar[dd]^{\pi}\ar[dr]^{\kappa} & 
        \\ 
        & (Z,\zeta)\ar[dl]^{\sigma}\\
        (Y,\eta)
        }        $$ commutes. Then $\sigma: (Z,\zeta)\longrightarrow (Y,\eta) $ is a relative $\mu$-boundary.    
\end{corollary}

\begin{example}
Let $\Gamma$ be a countable discrete group.
Let $\Sub(\Gamma)$ denote the collection of all subgroups of $\Gamma$ endowed with the Chabauty topology. Let $\lambda$ be a non-singular Borel probability measure on $\Sub(\Gamma)$. For instance, one can take an IRS (Invariant Random Subgroup). Let $\mu$ be a non-degenerate probability measure on $\Gamma$.
Consider
$$         \xymatrix{ \poi(\Gamma,\mu) \times (\Sub(\Gamma),\lambda)\ar[d]^{\kappa}
\\
        (Z,\zeta)\ar[d]^{\sigma}
        \\
        (\Sub(\Gamma),\lambda)
        }       $$
  such that $\sigma\circ\kappa=\pr_2$ is the projection to the second coordinate. 
 Then $\sigma:  (Z,\zeta) \longrightarrow (\Sub(\Gamma),\lambda) $ is a relative $\mu$-boundary by Corollary~\ref{cor: ift is rel prox}.
In particular, the so-called \textit{Bowen spaces}, introduced by \cite{Bowen} are relative $\mu$-boundaries over $(\Sub(\Gamma),\lambda)$.    
\end{example}

\subsection{Associated relative boundaries}
In general, we understand where to find instances of boundary behaviours. Furstenberg showed that measurable boundaries in the situation of stationary actions can be realised on the weak*-closure of the conditional measures (see \cite[Proposition~3.7]{HartKal}). Boundaries with topological characteristics come into existence when there's an affine action applied to a compact convex space (\cite[Theorem~III.2.3]{glasner1976proximal}). 

This suggests that if we want to associate a relative boundary map with a relatively stationary map, we need to look no further than the collection of probability measures (see \cite{GF09} for a non-relative version).
In particular, we prove the following.
\begin{theorem}\label{thm: bnd for X}
 Let  $\pi: (X,\nu)\longrightarrow (Y,\eta)$ be a relatively $\mu$-stationary $G$-factor map.
 Then there exists a Borel probability measure $\xi$ on $\prob(X)$ and a $\G$-factor map 
 $$\phi:(\prob(X), \xi)\longrightarrow (Y,\eta)$$ such that $\phi$ is a relative $\mu$-boundary. The measure $\xi$ is given by $$\xi:=\int_{Y}\int_{\G^{\mathbb{N}}} \delta_{\nu_{y,\omega}} \, d\mu^{\mathbb{N}}(\omega)\, d\eta(y), $$ and the factor map $\phi$ is defined by
     $$ \phi(\lambda):=y \text{ for } \supp(\lambda)\subseteq \pi^{-1}(y),
    $$ for $\lambda\in\prob(X)$ and $y\in Y$.
\end{theorem}

Before moving on to proving Theorem~\ref{thm: bnd for X}, we record the following simple observation for later use.

\begin{lemma}
\label{lem:compact}
Let $\pi:(X,\nu) \to (Y,\eta)$ be continuous surjective $G$-map. Then, the set
\begin{align*}
M = \bigcup_{y\in Y} \text{Prob}(\pi^{-1}(y))
\end{align*}
is compact.
\end{lemma}
\begin{proof}
It is enough to show that $M$ is closed. Towards this end, let $(\nu_k)_k\subset M$ be such that $\nu_k\to\nu$ in weak*-topology. We first observe that $\nu\in \text{Prob}(\pi^{-1}(y))$ if and only $\pi\nu=\delta_y$. Since $\nu_k\in M$, we see that $\pi\nu_k=y_k$ for some $y_k\in Y$. By passing to a subsequence if required, we can assume that $y_k\to y$. Then, using the continuity of $\pi$, we see that
\[\delta_y=\lim_k\delta_{y_k}=\lim_k\pi\nu_k=\pi\nu\]
Therefore, $\nu\in\text{Prob}(\pi^{-1}(y))$. The claim follows. 
\end{proof}

\begin{proof}[Proof of Theorem~\ref{thm: bnd for X}]
Let $\pi: (X,\nu)\longrightarrow (Y,\eta)$ be a relatively $\mu$-stationary $\G$-factor map. Without any loss of generality, by passing to the compact model, we assume that $X$ is compact and $\pi$ is continuous and surjective. Let $W=\prob(X)$ be the space of Borel probability measures on $X$. Let $$\xi:=\int_{Y}\int_{\G^{\mathbb{N}}} \delta_{\nu_{y,\omega}} \, d\mu^{\mathbb{N}}(\omega)\, d\eta(y), $$
where  $\delta_{\nu_{y,\omega}}(E)=1$ iff $\nu_{y,\omega}\in E$ and zero else.  
Note that $\bigcup_{y\in Y}\prob(\pi^{-1}(y))=W$ mod $\xi$. It follows from Lemma~\ref{lem:compact} that the set $\bigcup_{y\in Y}\prob(\pi^{-1}(y))$ is compact. Thus the limit $\nu_{y,\omega}$ lies inside this set and therefore
$$\xi\Big(\bigcup_{y\in Y}\prob(\pi^{-1}(y))\Big)=
\int_{Y}\int_{\G^{\mathbb{N}}} \delta_{\nu_{y,\omega}}\Big(\bigcup_{y\in Y}\prob(\pi^{-1}(y))\Big)
\, d\mu^{\mathbb{N}}(\omega)\, d\eta(y)=1.$$
Hence we can identify $(W,\xi)$ with $(\bigcup_{y\in Y} \prob(\pi^{-1}(y)),\xi)$. This enables us to define a map $\phi:  (W, \xi)\longrightarrow (Y,\eta)$ up to a null set by $\phi(w)=y$ if and only if $\supp(w)\subseteq \pi^{-1}(y)$. Note that this map is well-defined since the fibers $\{\pi^{-1}(y)\}_{y\in Y}$ are disjoint.

Let us first show that $\phi$ is relatively $\mu$-stationary. It is enough to establish that $$\int_{\G}\g\xi_{\g^{-1} y}\, d\mu(\g)=\xi_{y}$$ for $\eta$-a.e. $y\in Y$. 
    By Lemma~\ref{lem: equiv} we have 
    $\g\nu_{\g^{-1} y,\omega} =\nu_{y,\g \omega} $ for $\g\omega=\g(g_1,\ldots):=(\g,\g_1,\ldots)$. 
     From the definition of $\xi$, we obtain that $\xi_y=\int_{\G^{\mathbb{N}}} \delta_{\nu_{y,\omega}} \, d\mu^{\mathbb{N}}(\omega)$.
    Therefore,
\begin{align*}\int_{\G}\g\xi_{\g^{-1} y}\, d\mu(\g)
    &= \int_{\G}
    \int_{\G^{\mathbb{N}}} \g\delta_{\nu_{\g^{-1} y,\omega}} \, d\mu^{\mathbb{N}}(\omega)
    \, d\mu(\g)
    \\&=
    \int_{\G}
    \int_{\G^{\mathbb{N}}} \delta_{\g \nu_{\g^{-1} y,\omega}} \, d\mu^{\mathbb{N}}(\omega)
    \, d\mu(\g)
    \\&=\int_{\G}
    \int_{\G^{\mathbb{N}}} \delta_{ \nu_{ y,\g\omega}} \, d\mu^{\mathbb{N}}(\omega)
    \, d\mu(\g)    
     \\& =     \int_{\G^{\mathbb{N}}} \delta_{ \nu_{ y,\omega}} \, d\mu^{\mathbb{N}}(\omega)\\&= \xi_y.
    \end{align*}
It is left to show that $\phi$ is relatively proximal.
    Let  $\omega=(\g_1,\g_2,\ldots)\in\G^{\mathbb{N}}$ be such that the limit $\xi_{y,\omega}$ exists. Then, we obtain that
    \begin{align*}\xi_{y,\omega}&=\lim\limits_{n\to\infty} \g_1\ldots\g_n\xi_{\g_n^{-1}\ldots\g_1^{-1}y}
\\&=\lim\limits_{n\to\infty} \int_{\G^{\mathbb{N}}} \delta_{\g_1\ldots\g_n\nu_{\g_n^{-1}\ldots\g_1^{-1}y,\widetilde{\omega}}} \, d\mu^{\mathbb{N}}(\widetilde{\omega})
\\&=\lim\limits_{n\to\infty} \int_{\G^{\mathbb{N}}} \delta_{\nu_{y,\g_1\ldots\g_n\widetilde{\omega}}} \, d\mu^{\mathbb{N}}(\widetilde{\omega})&\text{(Lemma~\ref{lem: equiv})}\\&=\delta_{\nu_{y,\omega}}.
    \end{align*} Thus, $\phi$ is a relative boundary.
\end{proof}

\section{Relative Structure theorem}\label{sec: str}
The decomposition of dynamical systems into structured components is a powerful tool in ergodic theory. In this section, we prove a relative version of the Furstenberg-Glasner Structure Theorem, showing that any relative stationary system can be built from a relative boundary and a relatively measure-preserving map. In particular, we prove Theorem~\ref{thm: structure intro}. To this end, we collect some ingredients.
\subsection{Relative Joinings}\label{sec: join}
To formulate the structure theorem, we first need a mechanism to combine two relative dynamical systems over the same base. This motivates the definition of a relative joining.
For contractible spaces, relative joinings were considered by  
\cite{creutz2016contractive}.
\begin{definition}[Relative joining] \label{def: join}
Let $\pi_1: (X,\nu)\longrightarrow (Y,\eta)$ and $\pi_2: (Z,\zeta)\longrightarrow (Y,\eta)$ be relatively $\mu$-stationary $\G$-factor maps. We say that $\psi: (X\times Z,\iota)\longrightarrow (Y,\eta)$ is a \textit{relative joining} of $\pi_1$ and $\pi_2$, if 
$\psi=\pi_1\circ\pr_1=\pi_2\circ \pr_2$ and
$\psi$ is relatively $\mu$-stationary. 
\end{definition}

In the above definition, we obtain the following commutative diagram
     $$     \xymatrix{ &(X\times Z,\iota)\ar[dd]^{\psi}\ar[dl]^{\pr_1}\ar[dr]^{\pr_2} & \\
     (X,\nu)\ar[dr]^{\pi_1} & &(Z,\zeta)\ar[dl]^{\pi_2}
        \\ 
       & (Y,\eta)&
        }     ,   $$ where $\pr_i$ denotes the projection to the $i$-th coordinate.
Recall that we denoted \begin{equation}\label{eq: limagain}\nu_{y,\omega}:=\lim\limits_{n\to\infty} g_1\ldots g_n\nu_{g_n^{-1}\ldots g_1^{-1} y}    
\end{equation} 
for $\omega=(g_1,g_2,\ldots)\in\G^{\mathbb{N}}$ and $y \in Y$, whenever this limit exists in the weak*-topology on $\prob(X)$.
\begin{example}\label{ex: joining}
    Let $\pi_1:(X,\nu)\longrightarrow (Y,\eta)$
and $\pi_1:(Z,\zeta)\longrightarrow (Y,\eta)$ be continuous maps on compact Hausdorff spaces. Assume $\pi_1$ and $\pi_2$ are relatively $\mu$-stationary. 
Then the following defines a relative joining of $\pi_1$ and $\pi_2$.
Set $$\nu\vee \zeta:= \int_Y \int_{\G^{\mathbb{N}}} \nu_{y,\omega}\otimes \zeta_{y,\omega}\, d\mu^{\mathbb{N}}(\omega)\, d\eta(y)
$$ where the limit measures are defined by equation~\eqref{eq: limagain}.
Note that $\supp(\nu \vee \zeta)\subseteq \bigcup_{y\in Y} \pi_1^{-1}(y)\otimes\pi_2^{-1}(y),$ therefore we obtain $ \pi_1(\pr_1(x,z))=\pi_1(x)=\pi_2(z)=\pi_2( \pr_2(x,z))$
for $\nu\vee \zeta$-almost all $(x,z)\in X\otimes Z$.
Let us verify that $\psi:=\pi_1\circ\pr_1=\pi_2\circ \pr_2$ is a relatively $\mu$-stationary map from $(X\otimes Z,\nu\vee \zeta)$ to $(Y,\eta)$. Note that the action of $\G$ on $X\otimes Z$ is given by $\g(x,z)=(\g x,\g y)$. By construction, the disintegration measures are given by $(\nu\vee \zeta)_y= \int_{\G^{\mathbb{N}}} \nu_{y,\omega}\otimes \zeta_{y,\omega}\, d\mu^{\mathbb{N}}(\omega).$ Thus, using the equivariance Lemma~\ref{lem: equiv} of the limit measures, we obtain
\begin{align*}
\int_{\G} \g(\nu\vee \zeta)_{\g^{-1} y}\, d\mu(\g)
&=
\int_{\G} \int_{\G^{\mathbb{N}}} (\g\nu_{\g^{-1}y,\omega})\otimes (\g\zeta_{\g^{-1}y,\omega})\, d\mu^{\mathbb{N}}(\omega)  \, d\mu(\g)
\\&= 
\int_{\G} \int_{\G^{\mathbb{N}}} (\nu_{y,\g\omega})\otimes (\zeta_{y,\g \omega})\, d\mu^{\mathbb{N}}(\omega)  \, d\mu(\g)
\\&= 
 \int_{\G^{\mathbb{N}}} (\nu_{y,\omega})\otimes (\zeta_{y, \omega})\, d\mu^{\mathbb{N}}(\omega) 
\\&= (\nu\vee \zeta)_y,
\end{align*}
for $\eta$-a.e. $y\in Y$. This shows that $\psi$ is a relative joining in the sense of Definition~\ref{def: join}.

\end{example}
The above example has the following properties, similar to \cite[Proposition 3.1]{GF09} but in a relative version.
\begin{proposition}
\label{proposition:uniqejoining}
    The measure 
    \begin{equation}\label{eq: nu vee zeta}
        \nu\vee \zeta:= \int_Y \int_{\G^{\mathbb{N}}} \nu_{y,\omega}\otimes \zeta_{y,\omega}\, d\mu^{\mathbb{N}}(\omega)\, d\eta(y)
    \end{equation} of Example~\ref{ex: joining} has the following properties:
\begin{enumerate}
    \item\label{item: join prod} If $\pi_1$ is relatively measure-preserving, then $$\nu\vee \zeta=\int_{Y} \nu_y\otimes\zeta_y\, d\eta(y).$$
    \item\label{item: join uniq} If $\pi_1$ is relatively proximal, then $\nu\vee \eta$ is the only relative joining in the sense of Definition~\ref{def: join}.
\end{enumerate}
\end{proposition}

\begin{proof}
The proof of (\ref{item: join prod}) is straightforward: If $\pi_1$ is relatively measure-preserving then using Theorem~\ref{thm: lim ex}, we see that
\begin{align*}
\nu\vee \zeta&= \int_{Y} \int_{\G^{\mathbb{N}}} \nu_{y,\omega}\otimes \zeta_{y,\omega}\, d\mu^{\mathbb{N}}(\omega)\, d\eta(y)
\\&= 
\int_{Y} \nu_y\otimes\Big(\int_{\G^{\mathbb{N}}} \zeta_{y,\omega}\, d\mu^{\mathbb{N}}(\omega)\Big)\, d\eta(y)
\\&=
\int_{Y} \nu_y\otimes\zeta_y\, d\eta(y).
\end{align*}
To show (\ref{item: join uniq}), let $\pi_1$ be relatively proximal, i.e. $\nu_{y,\omega}=\delta_{\bnd(y,\omega)}$ for almost every $y\in Y$ and almost every $\omega\in \G^{\mathbb{N}}$.
Let $\lambda$ be any  relative joining of $\pi_1$ and $\pi_2$. We shall show that $\lambda=\nu \vee \zeta$ that is defined as above. Note that the limits 
$\lambda_{y,\omega}$
are probability measures on $X\times Y$. Moreover, by definition, $\pr_1\lambda=\nu$, thus $(\pr_1 \lambda)_y=\nu_y$ for almost every $y\in Y$. We claim that $(\pr_1 \lambda)_y=\pr_1(\lambda_y)$. Indeed,
$$ \int_Y \pr_1(\lambda_y)\, d\eta(y)= \pr_1\lambda=\nu$$ which proves the claim due to uniqueness of the disintegration. (Observe that the disintegrations are taken w.r.t. different maps, namely $\pi_1$ and $\psi=\pi\circ \pr_1$.)
Therefore,
$$\pr_1 (\lambda_{y,\omega})= (\pr_1\lambda)_{y,\omega}=\nu_{y,\omega}=\delta_{\bnd(y,\omega)}$$
for almost every $y\in Y$ and almost every $\omega\in\G^{\mathbb{N}}$.
Thus, 
$$\lambda_{y,\omega}=\delta_{\bnd(y,\omega)}\otimes \xi(y,\omega)$$ 
with some  probability measures $\xi(y,\omega)$ on $Z$ such that $$\zeta=\pr_2 \lambda=\int_Y \int_{G^{\mathbb{N}}}\xi(y,\omega)\, d\mu^{\mathbb{N}}(\omega) \, d\eta(y).$$
We claim that $\xi(y,\omega)=\zeta_{y,\omega}$. Indeed, as before, $$\xi(y,\omega)=\pr_2(\lambda_{y,\omega})=(\pr_2\lambda)_{y\omega}=\zeta_{y,\omega}.$$
Therefore, $$\lambda_y=\int_{\G^{\mathbb{N}}} \lambda_{y,\omega}\, d\mu^{\mathbb{N}}(\omega)=
\int_{\G^{\mathbb{N}}} \delta_{\bnd(y,\omega)} \otimes \zeta_{y,\omega}\, d\mu^{\mathbb{N}}(\omega)= 
\int_{\G^{\mathbb{N}}} \nu_{y,\omega} \otimes \zeta_{y,\omega}\, d\mu^{\mathbb{N}}(\omega),$$
which implies that $\lambda=\nu\vee \zeta$ defined as in (\ref{eq: nu vee zeta}), which was to prove.
\end{proof}
\subsection{Proof of the Structure Theorem}
Not every relative $\mu$-stationary space is a relative $\mu$-boundary. Thus, the question appears whether to every relative $\mu$-stationary map we can always associate a relative $\mu$-boundary map which is \say{not too far away}. In a non-relative version, this is Furstenberg-Glasner's Structure Theorem \cite[Theorem~4.3]{GF09}.

We start with the following definition of relating two-factor maps.

\begin{definition}[Relative 1-1 in the limit]\label{def: rel 1-1} 
Let $\pi_1: (X,\nu)\longrightarrow (Y,\eta)$ and $\pi_2: (Z,\zeta)\longrightarrow (Y,\eta)$ be relatively $\mu$-stationary $\G$-factor maps.
    A map $\phi:(X,\nu) \longrightarrow (Z,\zeta)$ is called  \textit{relatively 1-1 in the limit} if $$ \phi:(\pi_1^{-1}(y),\nu_{y,\omega}) \longrightarrow (\pi_2^{-1}(y),\zeta_{y,\omega}) \text{ is injective for $\mu^{\mathbb{N}}$-a.e. $\omega\in\G^{\mathbb{N}}$}
    $$ for $\eta$-a.e. $y\in Y$, with $\zeta_{y,\omega}=\phi \nu_{y,\omega} $
    and   
     $$     \xymatrix{ (X,\nu)\ar[dr]^{\pi_1}\ar[rr]^{\phi} & &(Z,\zeta)\ar[dl]^{\pi_2}
        \\ 
       & (Y,\eta)&
        }       $$ 
    commutes.
\end{definition}
In a non-relative version, i.e. for $(Y,\eta)=(\{\star\}, \delta_{\star})$, the above definition was introduced by Furstenberg-Glasner in \cite{GF09} under the name \textit{proximal map}. Since their name might cause confusion in our setting, we use a different one here.

We have now assembled all the necessary ingredients—relative boundaries, joinings, and the notion of maps that are 1-1 in the limit—to state and prove our main structural result. In particular, we prove a relative structure theorem, which generalises  \cite[Theorem~4.3]{GF09}.  (This is Theorem \ref{thm: structure intro} in the Introduction.)
\begin{theorem}\label{thm: structure}
    Let $\pi: (X,\nu)\longrightarrow (Y,\eta)$ be a relatively $\mu$-stationary $\G$-factor map. Then there exists a relative $\mu$-boundary map $\phi: (\prob(X),\xi)\longrightarrow (Y,\eta)$ such that 
$$     \xymatrix{ &(X\times \prob(X),\nu \vee \xi)\ar[dd]^{\psi}\ar[dl]^{\pr_1}\ar[dr]^{\pr_2} & \\
     (X,\nu)\ar[dr]^{\pi} & &(\prob(X),\xi)\ar[dl]^{\phi}
        \\ 
       & (Y,\eta)&
        }        $$
    commutes and $\pr_1$ is a relatively 1-1 in the limit and $\pr_2$ is relatively measure preserving.
\end{theorem}

\begin{proof} By Proposition~\ref{compactmodelsforfactormaps} we may assume that $X$ is compact Hausdorff and $\pi$ is continuous. 
The probability measure $\xi$ on $\prob(X)$ is given as in Theorem~\ref{thm: bnd for X}:
$$\xi:=\int_{Y}\int_{\G^{\mathbb{N}}} \delta_{\nu_{y,\omega}} \, d\mu^{\mathbb{N}}(\omega)\, d\eta(y), $$  where $\nu_{y,\omega}:=\lim\limits_{n\to\infty} g_1\ldots g_n\nu_{g_n^{-1}\ldots g_1^{-1} y} $ for a.e. $\omega=(\g_1,\g_2,\ldots)$. 
The map $\phi$ is given by $ \phi(w):=y \text{ for } \supp(w)\subseteq \pi^{-1}(y)$, which is almost everywhere well-defined, as shown in Theorem~\ref{thm: bnd for X}. Also, by  Theorem~\ref{thm: bnd for X}, we know that $\phi$ is a relative $\mu$-boundary. The measure on the product space $X\times \prob(X)$ is given by 
 $$\nu\vee \xi:= \int_Y \int_{\G^{\mathbb{N}}} \nu_{y,\omega}\times \xi_{y,\omega}\, d\mu^{\mathbb{N}}(\omega)\, d\eta(y).
$$ 
In Example~\ref{ex: joining} we have verified that $\pi\circ \pr_1=\phi\circ \pr_2=:\psi$  and $\psi$ is relatively $\mu$-stationary. It is left to show that the projections $\pr_1$ and $\pr_2$ are relatively 1-1 in the limit and relatively measure preserving, respectively. 
By definition of $\xi$ we see that $\xi_y=\int_{\G^{\mathbb{N}}} \delta_{\nu_{y,\omega}} \, d\mu^{\mathbb{N}}(\omega)$ and thus 
\begin{align*} \xi_{y,\omega}&=\lim\limits_{n\to\infty} \int_{\G^{\mathbb{N}}}\delta_{\g_1\ldots\g_n \nu_{\g_n^{-1}\ldots\g_1^{-1} y, \widetilde{\omega}}} \, d\mu^{\mathbb{N}}(\widetilde{\omega})\\&= 
\lim\limits_{n\to\infty} \int_{\G^{\mathbb{N}}}\delta_{ \nu_{ y, \g_1\ldots\g_n\widetilde{\omega}}} \, d\mu^{\mathbb{N}}(\widetilde{\omega})&\text{(using Lemma~\ref{lem: equiv})} 
\\&=\delta_{\nu_{y,\omega}}.
\end{align*} Therefore we  obtain \begin{equation}\label{eq: vee}
\nu\vee \xi= \int_Y \int_{\G^{\mathbb{N}}} \nu_{y,\omega}\times \delta_{\nu_{y,\omega}}\, d\mu^{\mathbb{N}}(\omega)\, d\eta(y).
\end{equation}
We claim that the disintegration with respect to $\pr_2$ is given by %\begin{equation}\label{eq: dis proj}
   $$ (\nu\vee \xi)_m=m\times \delta_{m} \text{ for $\xi$-a.e. }m\in\prob(X).$$
    %\end{equation} 
    Indeed, by the definition of $\xi$, we see that
$$
\int_{\prob(X)} m\times \delta_{m}\,  d\xi(m)=
\int_Y \int_{\G^{\mathbb{N}}} \nu_{y,\omega}\times \delta_{\nu_{y,\omega}}\, d\mu^{\mathbb{N}}(\omega)\, d\eta(y)=
\nu\vee \xi.
$$
Now, we can conclude that $$
 \g(\nu\vee \xi)_m=\g m\times \delta_{\g m} = (\nu\vee \xi)_{\g m},
$$ for all $\g \in \G$ and for $\xi$-a.e. $m\in\prob(X)$, hence $\pr_2$ is relatively measure preserving.
It is left to show that 
for $\eta$-a.e. $y\in Y$ and $\mu^{\mathbb{N}}$-a.e. $\omega\in\G^{\mathbb{N}}$ the projection to the first coordinate
$$\pr_1: (X\times \prob(X), (\nu\vee\xi)_{y,\omega}) \longrightarrow (X, \nu_{y,\omega})$$ is injective up to a null-set on $X\times \prob(X)$ with respect to $ (\nu\vee\xi)_{y,\omega}$.
By (\ref{eq: vee}) we see that $(\nu\vee\xi)_y=\int_{\G^{\mathbb{N}}} \nu_{y,\omega}\times \delta_{\nu_{y,\omega}}\, d\mu^{\mathbb{N}}(\omega)$ and thus similar to the calculations above, we obtain that 
$$
(\nu\vee\xi)_{y,\omega}=\nu_{y,\omega}\times \delta_{\nu_{y,\omega}}
$$
for $\eta$-a.e. $y\in Y$ and $\mu^{\mathbb{N}}$-a.e. $\omega\in\G^{\mathbb{N}}$. Note that this also implies that $ (\nu\vee\xi)_{y,\omega}\circ\pr_1^{-1}=\nu_{y,\omega}$, hence that $\pr_1$ is indeed a factor map with respect to these measures. Now if $\pr_1(x,m)=\pr_1(\widetilde{x},\widetilde{m})$ then $x=\widetilde{x}$ by definition of $\pr_1$. Since $(\nu\vee\xi)_{y,\omega}=\nu_{y,\omega}\times \delta_{\nu_{y,\omega}}$, we obtain that $m=\nu_{y,\omega}=\widetilde{m}$ up to a null set. This proves that $\pr_1$ is relatively 1-1 in the limit. This ends the proof.
\end{proof}

\section{The relative Poisson boundary}
\label{Sec: Poissonbdryrelative}
Among all relative boundaries, there is a unique maximal one, which we call the \textit{relative Poisson boundary}. We will make this precise in this section. In particular, we devote this section to the proof of Theorem~\ref{thm:Poissonrel_intro}.

Let $\poi(\G,\mu)$ denote the usual Poisson boundary of $(\G,\mu)$. We will show that the relative Poisson boundary over $(Y,\eta)$ is given by the projection of $\poi(\G,\mu)\otimes (Y,\eta)\longrightarrow (Y,\eta)$. 

\subsection{A construction}
We construct the relative Poisson boundary as shift ergodic components, similar to what is done for the usual Poisson boundary. Fix some standard $\G$-space $(Y,\eta)$. Let us consider the following shift map $$
T:(\G^{\mathbb{N}}\times  Y, \mu^{\mathbb{N}}\otimes \eta)\longrightarrow 
(\G^{\mathbb{N}}\times  Y, \mu^{\mathbb{N}}\otimes \eta)$$
$$T((\g_1,\g_2,\ldots),y):=((\g_1 \g_2,\g_3,\ldots), y)
$$ which acts as the usual shift on the sequence space $\Omega$ and leaves $Y$ untouched. Now, $T$ commutes with the diagonal action of $\G$ on $\G^{\mathbb{N}}\times  Y$  given by $$\g\cdot ((\g_1,\g_2,\ldots),y):=((\g\g_1,\g_2,\ldots),\g y).$$
Indeed,
$$ T(\g\cdot ((\g_1,\g_2,\ldots),y))=  
((\g\g_1\g_2,\g_3\ldots),\g y)=
\g(T((\g_1,\g_2,\ldots),\g y)).
$$
Therefore, we may take the $T$-ergodic components and still obtain a $\G$-action on this new space. The $T$-ergodic components are given by the Mackey point realisation with respect to the sub-$\sigma$-algebra of $T$-invariant measurable sets. Since $T$ does not change the $Y$-component, we obtain that this space is of the form $$\poi(\G,\mu)\otimes (Y,\eta),$$ where $\poi(\G,\mu)$ denotes the usual Poisson boundary of the random walk $\mu$ on $\G$.

Let us now verify that $\pr_2:\poi(\G,\mu)\otimes (Y,\eta)\longrightarrow (Y,\eta) $ is a relative $\mu$-boundary in the sense of Definition~\ref{def:  rel bnd}. Let us denote $\poi(\G,\mu)=:(X,\nu)$.
The disintegration measures of $\pr_2$ are given by  $(\nu\otimes\eta)_y=\nu\otimes \delta_{y}
$ for $\eta$-a.e. $y\in Y$. Thus, relative $\mu$-stationarity follows by
$$\int_{\G} \g (\nu\otimes\eta)_{\g^{-1} y}\, d\mu(\g)=
\int_{\G} (\g \nu)\otimes(\g \delta_{\g^{-1}y})\, d\mu(\g)
$$
$$=\int_{\G} (\g \nu)\otimes \delta_{y}\, d\mu(\g)=\nu\otimes \delta_{y}$$
due to $\mu$-stationarity of $\nu$.
To see that $\pr_2$ is relative $\mu$-proximal, we consider
$$(\nu\otimes\eta)_{y,\omega}=\lim\limits_{n\to\infty}\g_1\cdots \g_n(\nu\otimes \delta_{\g_n^{-1}\ldots\g_1^{-1}y})=\nu_{\omega}\otimes \delta_y= \delta_{(\bnd(\omega),y)}
$$ where we used that  $\nu_{\omega}:=\lim\limits_{n\to\infty} \g_1\cdots \g_n \nu=\delta_{\bnd(\omega)}$ for some element $\bnd(\omega)\in X$ for $\mu^{\mathbb{N}}$-a.e. $\omega=(\g_1,\g_2,\ldots)\in\G^{\mathbb{N}}$ by the proximality properties of the usual Poisson boundary.
This shows that  $\pr_2:\poi(\G,\mu)\otimes (Y,\eta)\longrightarrow (Y,\eta) $ is a relative $\mu$-boundary.

\subsection{Maximality and uniqueness}
We have constructed a candidate for the relative Poisson boundary using the shift space. We now prove that this object satisfies the universal property of maximality, meaning any other relative boundary is a factor of it.
\begin{proposition}[Maximality]\label{proposition: max poi}
Let  $\phi: (Z,\zeta)\longrightarrow (Y,\eta)$ be any  relative $\mu$-boundary. Then there exists a $\G$-factor map $\pi:\poi(\G,\mu)\otimes (Y,\eta)\longrightarrow (Z,\zeta) $  such that $\pr_2=\phi\circ\pi $, i.e. the diagram
$$     \xymatrix{ \poi(\G,\mu)\otimes (Y,\eta)\ar[dd]^{\pr_2}\ar@{..>}[dr]^{\pi} & 
\\
   &  (Z,\zeta)\ar[dl]^{\phi} 
        \\ 
       (Y,\eta)&        }        $$
commutes.  
\end{proposition}

\begin{proof} We fist pass to compact models and assume that $\phi$ is continuous by Proposition~\ref{compactmodelsforfactormaps}.
Consider the almost everywhere defined map $\bnd: \G^{\mathbb{N}}\times Y\longrightarrow Z $ defined via the relative $\mu$-proximality $\zeta_{\omega,y}=\delta_{\bnd(\omega,y)}$ for $\mu^{\mathbb{N}}\otimes \eta$-a.e. $(\omega,y)\in \G^{\mathbb{N}}\times Y$. Let us denote $\poi(\G,\mu)=(X,\nu)$.
Since $\bnd:\G^{\mathbb{N}}\times Y\longrightarrow Z $ is $\G$-equivariant and measurable (see Corollary~\ref{corollary: bnd map}),
we obtain a $\G$-equivariant measurable map $$\pi:\poi(\G,\mu)\otimes (Y,\eta)\longrightarrow (Z, \pi(\nu\otimes\eta)),$$ when taking the $T$-ergodic components, i.e. $\pi([\omega]_{\sim},y):=\bnd(\omega,y)$ for the equivalence class $[\omega]_{\sim}$ defined by Mackey's point-realisation of the $T$-invariant components.
Moreover, $\supp(\zeta_{\omega,y})\subseteq \phi^{-1}(y)$, thus  $\bnd(\omega,y)\in \phi^{-1}(y)$ and therefore $\phi(\bnd(\omega,y))=y$ which implies that $\phi\circ \pi= \pr_2$.  Therefore, $\phi\circ\pi$ is a surjective map from $X\times Y$ to $Y$ up to null sets.
It is left to show that $\pi(\nu\otimes\eta)=\zeta$.
Recall that $\zeta_y=\int_{\G^{\mathbb{N}}}\zeta_{\omega,y}\, d\mu^{\otimes\mathbb{N}}(\omega)$ for $\eta$-a.e. $y\in Y$ by the barycenter equation in Theorem~\ref{thm: lim ex}. Thus,
$$\zeta= \int_Y \zeta_y\, d\eta(y)=
\int_Y \int_{\G^{\otimes\mathbb{N}}}\delta_{\bnd(\omega,y)} \, d\mu^{\otimes\mathbb{N}}(\omega)\, d\eta(y) =
(\mu^{\otimes\mathbb{N}}\otimes  \eta)\circ \bnd^{-1}.
$$
Now, since $\bnd$ is $\G$-equivariant, we obtain for every $A\in\mathcal{B}(Y)$ that
$$
(\nu\otimes \eta)(\pi^{-1}(A))
$$ $$=(\nu\otimes \eta)(\{([\omega]_{\sim},y)\in X\times Y\, :\, \bnd(\omega,y)\in A\})
$$
$$
= (\mu^{\otimes\mathbb{N}}\otimes \eta)(\{(\omega,y)\in \G^{\mathbb{N}}\times Y\, :\, \bnd(\omega,y)\in A\})
$$
$$= (\mu^{\mathbb{N}}\otimes \eta) (\bnd^{-1}(A)),
$$
hence $(\mu^{\mathbb{N}}\otimes \eta )\circ \bnd^{-1}=(\nu\otimes \eta)\circ \pi^{-1} $ which finishes the proof.
\end{proof}
The maximality property, combined with the rigidity of proximal maps, allows us to deduce the uniqueness of the relative Poisson boundary. In particular, applying  Proposition~\ref{proposition: id map} to the relative Poisson boundary together with the maximality of Proposition~\ref{proposition: max poi}, we can conclude the uniqueness of the relative Poisson boundary up to $\ G$-isomorphisms.
\begin{corollary}[Uniqueness]\label{corollary: unique poi}
    The relative Poisson boundary is unique up to $\G$-isomorphism. Meaning that if there is any other "maximal"  $\mu$-boundary $\phi: (Z,\zeta)\longrightarrow (Y,\eta)$ in the sense of Proposition~\ref{proposition: max poi}, then there exists a $\G$-isomorphism $\psi: \poi(\G,\mu)\otimes (Y,\eta) \longrightarrow (Z,\zeta)$ 
    such that $\phi \circ \psi= \pr_2.$    
\end{corollary}

\begin{proof}
Assume $\phi: (Z,\zeta)\longrightarrow (Y,\eta)$ is a relative $\mu$-boundary such that every other boundary is a factor of it as in Proposition~\ref{proposition: max poi}. Then, in particular there exists a $\G$-factor map $\pi:(Z,\zeta)\longrightarrow \poi(\G,\mu)\otimes (Y,\eta)$ such that $\pr_2\circ \pi=\phi$.
On the other hand by Proposition~\ref{proposition: max poi}, there exists $\rho:\poi(\G,\mu)\otimes (Y,\eta) \longrightarrow (Z,\zeta)$ such that $\phi\circ \rho=\pr_2.$
    Thus  $$\alpha:=\pi\circ\rho: \poi(\G,\mu)\otimes (Y,\eta) \longrightarrow \poi(\G,\mu)\otimes (Y,\eta)$$
is a $\G$-equivariant factor map with $\alpha (\nu\otimes \eta)=(\nu\otimes \eta)$, 
where we denote $(X,\nu)=\poi(\G,\mu)$ as above. 
With the above we have $\pr_2\circ \alpha=\pr_2\circ \pi\circ \rho= \phi\circ \rho=\pr_2$. Thus, we can apply  
 Proposition~\ref{proposition: id map} and obtain  $\alpha=id$, which ends the proof.
\end{proof}

%\subsection{Properties}
\begin{remark}
    The action of $\G$ on the relative Poisson boundary is amenable because the action of $\G$ on $\poi(\G,\nu)$ is. However, the action on the relative Poisson boundary might not be ergodic.
\end{remark}
Finally, we provide an example of a relative boundary that is not a product space, illustrating the non-trivial nature of these extensions.
\begin{example}
Let $(Y,\eta)$ be a $G$-space. Let $J\triangleleft L^{\infty}(Y,\eta)$ be a non-trivial $G$-invariant WOT-closed ideal. Consider 
\[\mathcal{M}=L^{\infty}(Y,\eta)+J\overline{\otimes} L^{\infty}(\text{Poi}(G,\mu),\nu_B).\]
Combining Proposition~\ref{proposition: IF} along with Proposition~\ref{proposition: IF prox}, we see that $\mathcal{M}$ is a relative $\mu$-boundary over $(Y,\eta)$. We claim that $\mathcal{M}$ does not split. 

Let $\mathbb{E}_{\nu_B}:L^{\infty}(Y,\eta)\overline{\otimes}L^{\infty}(\text{Poi}(G,\mu),\nu_B)\to L^{\infty}(Y,\eta)$ be the conditional expectation defined by 
\[\mathbb{E}_{\nu_B}(f\otimes g)=\nu_B(f)g,~g\in L^{\infty}(\text{Poi}(G,\mu),\nu_B),\  f\in  L^{\infty}(Y,\eta).\]
Similarly,
$\mathbb{E}_{\eta}:L^{\infty}(Y,\eta)\overline{\otimes}L^{\infty}(\text{Poi}(G,\mu),\nu_B)\to L^{\infty}(\text{Poi}(G,\mu),\nu_B)$ is the conditional expectation defined by 
\[\mathbb{E}_{\eta}(f\otimes g)=f\eta(g),\ g\in L^{\infty}(\text{Poi}(G,\mu),\nu_B),\  f\in  L^{\infty}(Y,\eta).\]
We also observe that for $\sum_{i=1}^n j_i\otimes f_i\in J\otimes L^{\infty}(\text{Poi}(G,\mu),\nu_B)$, $$\mathbb{E}_{\nu_B}\left(\sum_{i=1}^n j_i\otimes f_i\right)=\sum_{i=1}^nj_i \nu_B(f_i)\in J.$$
Since $\mathbb{E}_{\nu_B}$ is normal, it follows that $\mathbb{E}_{\nu_B}(b)\in J$ for $b\in J\otimes L^{\infty}(\text{Poi}(G,\mu),\nu_B)$.

To establish that $\mathcal{M}$ does not split, we shall prove that
\begin{enumerate}
\item\label{emptyintersection} $\mathcal{M}\cap L^{\infty}(\text{Poi}(G,\mu),\nu_B)=\mathbb{C}$
    \item \label{everything}$\mathbb{E}_{\eta}\left(\mathcal{M}\right)=L^{\infty}(\text{Poi}(G,\mu),\nu_B)$.
\end{enumerate}
Let us first prove (\ref{emptyintersection}). We begin by observing that $$\left(J\overline{\otimes}L^{\infty}(\text{Poi}(G,\mu),\nu_B)\right)\cap L^{\infty}(\text{Poi}(G,\mu),\nu_B)=\{0\}.$$
Indeed, let $f\in L^{\infty}(\text{Poi}(G,\mu),\nu_B)$ be such that it is in $J\overline{\otimes}L^{\infty}(\text{Poi}(G,\mu),\nu_B)$. Since $$\mathbb{E}_{\nu_B}\left(J\overline{\otimes}L^{\infty}(\text{Poi}(G,\mu),\nu_B)\right)\subset J,$$ it follows that $$\nu_B(f^*f)=\mathbb{E}_{\nu_B}(f^*f)\in J.$$
If $\nu_B(f^*f)\ne 0$, then we would obtain that $J=L^{\infty}(Y,\eta)$ which would contradict the non-triviality of $J$. Therefore, $\nu_B(f^*f)=0$ which shows that $f=0$ $\nu_B$-almost everywhere.

Now, let $f\in \mathcal{M}\cap L^{\infty}(\text{Poi}(G,\mu),\nu_B)$. Then, $f=a+b$ for some $a\in L^{\infty}(Y,\eta)$ and $b\in J\overline{\otimes} L^{\infty}(\text{Poi}(G,\mu),\nu_B)$.  Applying the conditional expectation $\mathbb{E}_{\nu_B}$ to $f=a+b$, we see that 
\[\nu_B(f)=a+\mathbb{E}_{\nu_B}(b).\]
Therefore,
\[f-\nu_B(f)=b-\mathbb{E}_{\nu_B}(b).\]
Since $b-\mathbb{E}_{\nu_B}(b)\in J\overline{\otimes}L^{\infty}(\text{Poi}(G,\mu),\nu_B)$, we see that $f-\nu_B(f)\in J\overline{\otimes}L^{\infty}(\text{Poi}(G,\mu),\nu_B)$. However, $f-\nu_B(f)\in L^{\infty}(\text{Poi}(G,\mu),\nu_B)$. Therefore, from our first observation, we see that $f-\nu_B(f)=0$. This proves (\ref{emptyintersection}).

Let us now prove (\ref{everything}). Observe that for an element of the form $j\otimes f\in J\overline{\otimes} L^{\infty}(\text{Poi}(G,\mu),\nu_B)$ with $\eta(j)\ne 0$, $$\mathbb{E}_{\eta}\left(\frac{j}{\eta(j)}\otimes f\right)=1\otimes f.$$
This establishes the claim. 
\end{example}
The construction above has been used in \cite{zacharias2001splitting} to construct intermediate algebras which do not split.

\end{document}